\documentclass[pdflatex,sn-mathphys-num]{sn-jnl}

\usepackage{graphicx}%
\usepackage{multirow}%
\usepackage{amsmath,amssymb,amsfonts}%
\usepackage{amsthm}%
\usepackage{mathrsfs}%
\usepackage[title]{appendix}%
\usepackage{xcolor}%
\usepackage{textcomp}%
\usepackage{manyfoot}%
\usepackage{booktabs}%
\usepackage{algorithm}%
\usepackage{algorithmicx}%
\usepackage{algpseudocode}%
\usepackage{listings}%

\theoremstyle{thmstyleone}%
\theoremstyle{thmstyletwo}%

\theoremstyle{thmstylethree}%

\usepackage{subcaption}
\newtheorem{lemma}{Lemma}
\usepackage{csquotes}
\usepackage{bm}
\counterwithin{figure}{section}
\counterwithin{table}{section}
\def\div{\mbox{{\rm div}}}
\newcommand{\bA}{\mathbf{A}}
\newcommand{\bB}{\mathbf{B}}
\newcommand{\bC}{\mathbf{C}}
\newcommand{\bb}{\bm{b}}
\newcommand{\bd}{\bm{d}}
\newcommand{\bx}{\bm{x}}

\newcommand{\bo}{\bm{0}}
\newcommand{\B}{\mathcal{B}}
\newcommand{\F}{\mathcal{F}}
\newcommand{\G}{\mathcal{G}}

\newcommand{\R}{\mathbb{R}}
\newcommand{\ddx}{\,dx}
\newcommand{\ddt}{\,dt}
\newcommand{\dx}{\Delta x}
\newcommand{\dy}{\Delta y}
\newcommand{\dt}{\Delta t}
\newcommand{\Id}{\mathbf{Id}}
\newcommand{\be}{\begin{equation}}
\newcommand{\ee}{\end{equation}}
\newcommand{\ba}{\begin{array}}
\newcommand{\ea}{\end{array}}
\newcommand{\prox}{\ensuremath{\operatorname{\mathbf{prox}}}}

\usepackage{glossaries}
\glsdisablehyper 
\newacronym{cg}{CG}{conjugate gradient}
\newacronym{pde}{PDE}{partial differential equation}
\newacronym{amg}{AMG}{algebraic multigrid}
\newacronym{mfg}{MFG}{mean field game}
\newacronym{dst}{DST}{discrete sine transform}
\newacronym{illt}{iLLT}{incomplete Cholesky factorization}
\newacronym{cp}{CP}{Chambolle--Pock}
\newacronym{fft}{FFT}{fast Fourier transform}
\newacronym{dct}{DCT}{discrete cosine transform}
\newacronym{dft}{DFT}{discrete Fourier transform}

\begin{document}

\title[Parallel-in-Time Preconditioning for Time-Dependent Variational Mean Field Games]{Parallel-in-Time Preconditioning for Time-Dependent Variational Mean Field Games}


\author*[1]{\fnm{Heidi Wolles} \sur{Lj\'{o}sheim}}\email{h.w.ljosheim@sms.ed.ac.uk}

\author[2]{\fnm{Dante} \sur{Kalise}}\email{d.kalise-balza@imperial.ac.uk}
\equalcont{These authors contributed equally to this work.}

\author[1]{\fnm{John W.} \sur{Pearson}}\email{j.pearson@ed.ac.uk}
\equalcont{These authors contributed equally to this work.}

\author[3]{\fnm{Francisco J.} \sur{Silva}}\email{francisco.silva@unilim.fr}
\equalcont{These authors contributed equally to this work.}

\affil*[1]{\orgdiv{School of Mathematics and Maxwell Institute for Mathematical Sciences}, \orgname{The University of Edinburgh}, \orgaddress{\street{James Clerk Maxwell Building, The King's Buildings}, \city{Edinburgh}, \postcode{EH9 3FD}, \country{United Kingdom}}}

\affil[2]{\orgdiv{Department of Mathematics}, \orgname{Imperial College London}, \orgaddress{\street{South Kensington Campus}, \city{London}, \postcode{SW7 2AZ}, \country{United Kingdom}}}

\affil[3]{\orgdiv{Institut de recherche XLIM-DMI, UMR 7252 CNRS, Facult\'e des Sciences et Techniques}, \orgname{Universit\'e de Limoges}, \orgaddress{\city{Limoges}, \postcode{87060}, \country{France}}}


\abstract{We study the numerical approximation of a time-dependent variational mean field game system with local couplings and either periodic or Neumann boundary conditions. Following a variational approach, we employ a finite difference discretization and solve the resulting finite-dimensional optimization problem using the Chambolle--Pock primal--dual algorithm. As this involves computing proximal operators and solving ill-conditioned linear systems at each iteration, we embed within our solver a general class of parallel-in-time preconditioners based on suitably-chosen diagonalization techniques, applied using discrete Fourier transforms. These enable efficient, scalable iterative solvers for each linear system, with robustness across a wide range of viscosities. We further develop fast solvers for the resulting ill-conditioned systems arising at each time step, using exact recursive schemes for structured grids while allowing for other geometries. Numerical experiments confirm the improved performance and parallel scalability of our approach.}




\maketitle

\section{Introduction}
\Gls{mfg} theory has been introduced independently by J.-M. Lasry and P.-L. Lions~\cite{LasryLions06i,LasryLions06ii,LasryLions07} and by M. Huang, R. Malham\'e, and P. Caines~\cite{huang_large_2006}. The aim of \glspl{mfg} is to model the asymptotic behavior of equilibria of symmetric differential games as the number of players tends to infinity.
We refer the reader to~\cite{Gomes_survey,Carmona_Delarue_i,Carmona_Delarue_ii,Achdou_et_al_MFG_CIME} for background in \gls{mfg} theory, the numerical approximation of equilibria, probabilistic techniques, and applications. 

In this work we consider the following time-dependent \gls{mfg} system with local couplings and periodic boundary conditions:
\begin{equation}
\label{MFGcont_torus}
\left\{
\begin{aligned}
&-\partial_t u -\nu\Delta u + H(x, \nabla u) =f(x,m(x,t)) &\text{ in }\mathbb{T}^{d}\times [0,T]\,,\\
&\partial_t m -\nu\Delta m - \mbox{div}(\nabla_p H(x,\nabla u)m)=0 &\text{ in } \mathbb{T}^{d}\times [0,T]\,,\\
&m(x,0)=m_0(x), \hspace{0.3cm} u(x,T)=g(x,m(x,T)) &\text{ in }\mathbb{T}^{d}\,,
\end{aligned}
\right.\tag{${\rm MFG}_{\mathbb{T}^{d}}$}
\end{equation}
where $\mathbb{T}^{d}$ is the $d$-dimensional torus, $T>0$, the viscosity $\nu\geq 0$, $H:\mathbb{T}^{d}\times\R^d\to\R$ is jointly continuous, convex, and differentiable with respect to its second variable, $f$, $g: \mathbb{T}^{d}\times[0,+\infty)\to \R$ are continuous functions, and $m_0 \in L^1(\mathbb{T}^{d})$ satisfies $m_0 \geq 0$ a.e. and $\int_{\mathbb{T}^{d}}m_0(x) \ddx=1$, hence representing a probability density function. 
Given an open bounded domain $\emptyset\neq\Omega\subseteq\R^{d}$, we also consider a variant of $({\rm MFG}_{\mathbb{T}^{d}})$ with the following additional Neumann boundary conditions:
\begin{align*}
\nabla u(x)\cdot n(x) &= 0 \;\text{ on }\; \partial\Omega\times [0,T], \\
\left(\nu\nabla m+\nabla_p H(x,\nabla u)m\right)\cdot n(x)&=0 \;\text{ on } \;\partial\Omega\times [0,T],
\end{align*}
where now $H:\overline{\Omega}\times\R^d\to\R$, and $n(x)$ denotes a unit outward normal vector at $x\in\partial\Omega$. We also define $f$, $g: \overline{\Omega}\times[0,+\infty)\to \R$, $m_0:\R^{d}\to\R$, with support contained in $\overline{\Omega}$, and we replace ${\mathbb{T}^d}$ by $\Omega$ wherever it appears in the original system \eqref{MFGcont_torus}, denoting this new system by $({\rm MFG}_{\Omega})$.

Under various assumptions on the Hamiltonian $H$, the coupling functions $f$ and $g$, and the initial distribution $m_0$, the existence and uniqueness of solutions to~\eqref{MFGcont_torus} have been studied in several works, starting from~\cite{LasryLions06ii,LasryLions07}, using fixed point techniques. Instead of providing an exhaustive list of references, we refer the reader to~\cite{Gomes_regularity_book,MR4214774} and the references therein for a detailed account on the well-posedness of~\eqref{MFGcont_torus}. Regarding system~(MFG$_{\Omega})$, we refer the reader to~\cite{Cirant_Neumann,Porretta_ARMA} for existence and uniqueness results and to~\cite{Ricciardi_master_equation} for the study of the associated master equation.

The numerical approximation of \gls{mfg} systems has been an active research topic during the last decade. We refer the reader to~\cite{Achdou_Lauriere_review, osborne_analysis_2024, osborne_finite_2025} and the references therein for an overview of this subject. With classical sources employing the Newton method \cite{MR2679575,Achdou_Lauriere_review}, recent numerical solution methods for \glspl{mfg} include an augmented Lagrangian method \cite{MR3395203,andreev_17}, a primal--dual method \cite{BAKS,briceno-arias_implementation_2019}, and fictitious play, or, the generalized conditional gradient method \cite{MR4030259,Lavigne_Pfeiffer,MR4886523}. As noticed in~\cite{LasryLions07}, systems of the form~\eqref{MFGcont_torus}--(MFG$_{\Omega})$ admit a  variational formulation, i.e., they correspond to the optimality system of an underlying optimization problem. Taking advantage of this structure, a promising approach in numerical methods for \glspl{mfg} is to discretize the underlying \gls{pde} and treat it as the optimality system of a finite-dimensional problem, allowing the use of convex optimization techniques to calculate equilibria (see, e.g., \cite{MR2888257,MR3395203,Almulla_Ferreira_Gomes,BAKS,briceno-arias_implementation_2019,Liu_Osher_Nurbekyan,Lavigne_Pfeiffer,Briceno_Deride_Lopez_Silva,Briceno_Silva_Yang}).

As in~\cite{BAKS,briceno-arias_implementation_2019}, we focus on the application of the primal--dual method introduced in~\cite{MR2782122} applied to a finite difference discretization of \eqref{MFGcont_torus}. We also consider a similar discretization of the underlying variational problem associated with~(MFG$_{\Omega})$. At each iteration, the primal--dual method requires the computation of two so-called proximal operators \cite{Moreau_i,Moreau_ii}. Since these operators can be efficiently calculated when the Hamiltonian is given by
\begin{equation}
\mathbb{T}^{d}\times\R^{d}\ni(x,p)\mapsto H(x,p)=\frac{1}{q'}|p|^{q'}\in\R,
\label{power_hamiltonian}
\end{equation}
where $q'\in(1,+\infty)$ and $|p|$ denotes the Euclidean norm of $p$, in this work we restrict our attention to this type of Hamiltonian and we concentrate our efforts in providing an efficient solution method for the linear systems appearing in the primal--dual iterations.

The purpose of this paper is to improve the performance of the solution method in \cite{briceno-arias_implementation_2019} by considering a novel set of preconditioners for the linear systems arising in the primal--dual algorithm. 
Let $\mathcal{D}$ denote the discrete version of the differential operator $(\partial_t - \nu\Delta, \text{div})$, then the dual proximal operator inevitably involves applying the inverse of a matrix of the form $\mathcal{D}\mathcal{D}^*$.
Since non-zero viscosity \glspl{mfg} admit a second-order differential operator, this leads to a fourth-order operator in the matrix $\mathcal{D}\mathcal{D}^*$. The condition number for a similar (simpler) fourth-order operator was observed to scale as $\mathcal{O}(h^{-4})$ when the discretization parameter $h$ tends to zero in \cite{AchdouPerez2012}, showing the necessity of preconditioning in order to obtain accurate results in a reasonable time.

A number of techniques have been tested to improve the convergence of Krylov subspace solvers for systems of the form $\mathcal{D}\mathcal{D}^*$ arising from \glspl{mfg}, most popularly multigrid methods in \cite{AchdouPerez2012, andreev_17, briceno-arias_implementation_2019}. These show reasonable results which generally improve when increasing the diffusion $\nu$. More elaborate preconditioning strategies include applying a temporal transformation that produces a block-diagonalizable preconditioner \cite{andreev_17}. However, we note that \cite{andreev_17} deals with operators from a temporal discretization similar to continuous Galerkin time-stepping, and the alternating direction method of multipliers (ADMM) is employed as the solution method.

In this work, we propose a class of fully parallel-in-time preconditioners based on the theory of (block-) diagonalizing in time using \glspl{dft}, for systems arising when solving variational mean field games using the Chambolle--Pock algorithm. This allows a parallel implementation, which is currently the main advantage of the primal--dual algorithm. Moreover, the preconditioner is exact for zero viscosity, as observed in \cite{MR3158785}. In comparison, a second-order method like the Newton method will generally not permit a parallel-in-time preconditioner unless the equations are approximated in a way that yields time-uniform blocks. Indeed, a Newton linearization will lead to very different linear systems at every iteration, which is not the case with the Chambolle--Pock algorithm, meaning that we may re-use a preconditioner from one Chambolle--Pock iteration to another. Moreover, the Newton method can be shown to struggle with small viscosity regimes as $m$ may be negative at some iterations \cite{AchdouPerez2012}. This behaviour may be remedied by a method such as warm-starting the algorithm with smaller and smaller viscosity initial solutions or by using an intermediate solution found with a first-order method.

The preconditioning strategy is complete with a number of different options for solving the system involving the fourth-order differential operator which arises at each time step: we construct exact recursive linear algebra solvers for the finite difference schemes with either periodic or Neumann boundary conditions, and have the flexibility to defer to LU factorization and \gls{amg} methods for more general meshes.
In a number of numerical experiments, the preconditioners show robustness with respect to the viscosity parameter $\nu$, and good scaling with respect to their parallelization capabilities.

The key contributions of this paper are summarised as follows:
\begin{itemize}
    \item A framework for the preconditioned iterative solution of systems arising from variational \glspl{mfg}, that allows the application of different parallelization strategies across the time variable. Between these strategies, rapid convergence is achieved in settings with both large and small viscosity parameters.
    \item The incorporation of effective strategies for preconditioning the matrices across the spatial variables arising at every time-step, including a recursive, block \gls{dft} linear solver for certain settings.
    \item The practical and robust numerical solution of very large \gls{mfg} systems, previously untested in the literature to our knowledge, in particular systems within the numerically challenging regime of small but non-zero viscosities. 
    \item A fully parallel implementation for different \glspl{mfg}, leading to convincing strong and weak scaling results for our preconditioned iterative solver.
    \item The preconditioners are applicable beyond the Chambolle--Pock algorithm as systems involving the matrix $\mathcal{D}\mathcal{D}^*$ may arise in a number of solution algorithms, including when solving normal equations.
\end{itemize}

This paper is structured as follows: In Section \ref{sec:discrete}, we introduce some standard notation and recall the finite difference scheme for \eqref{MFGcont_torus} as originally proposed in \cite{MR2679575} and extend this to the variational formulation that we examine, including the setting with Neumann boundary conditions. Section \ref{sec:alg} is devoted to the primal--dual algorithm introduced in \cite{MR2782122}, and we discuss its application to the variational problem. In Section \ref{sec:prec}, we present the new preconditioning strategy. In Section \ref{sec:num}, we present numerical results for two well-known test problems, accompanied by a number of parallel scaling tests.

\section{Finite difference discretizations} \label{sec:discrete}
In this section, we first recall the finite difference discretization of~\eqref{MFGcont_torus} introduced in~\cite{MR2679575} and provide a slight variant to deal with the Neumann boundary conditions in~(MFG$_{\Omega})$. For the sake of simplicity, we will focus on the two-dimensional case $d=2$ and, as mentioned in the introduction, we consider Hamiltonians of the form~\eqref{power_hamiltonian}.

Consider a rectangular domain $\Omega=[a,b]\times[c,d]$, the height and width of which for periodic boundary conditions correspond to the circumferences of the torus. Given positive integers $N_t$, $N_x$, $N_y$ and setting $\dt=T/N_t$, $h=(\Delta x, \Delta y)=((b-a)/N_x,(d-c)/N_y)$, we define the time and space grids, respectively,
\begin{align*}
    T_{\dt}&:= \{t_k = k\dt \; ; \; k=0,\dots,N_t\}, \\
    \mathbb{T}_{h}^{2}&:= \{x_{i,j}=(a+i\dx,c+j\dy) \; ; \; i=0,1,\ldots,N_x-1,\;j=0,1,\ldots,N_y-1\}.
\end{align*}
For a discrete scalar field $\alpha:\mathbb{T}_{h}^{2}\times T_{\dt} \to\R$, indexed as $\alpha_{i,j}^k$, we define for $T_{\dt}$ and $\mathbb{T}_{h}^{2}$ the finite difference operators
\begin{equation*}
\begin{aligned}
&(D_1 \alpha)^k_{i,j} := \frac{\alpha^k_{i+1,j}-\alpha^k_{i,j}}{\dx},\quad (D_2 \alpha)^k_{i,j} := \frac{\alpha^k_{i,j+1}-\alpha^k_{i,j}}{\dy},\quad 
(D_t \alpha)^k_{i,j} := \frac{\alpha^{k+1}_{i,j}-\alpha^k_{i,j}}{\dt}\,,\\
&\text{and } (\Delta_h \alpha)^k_{i,j} := \frac{\alpha_{i+1,j}^{k}-2\alpha^k_{i,j}+\alpha_{i-1,j}^{k}}{(\dx)^2}+ \frac{\alpha_{i,j+1}^{k}-2\alpha^k_{i,j}+\alpha_{i,j-1}^{k}}{(\dy)^2},
\end{aligned}
\end{equation*}
where for periodic boundary conditions the algebraic operations over the indices $i$ and $j$ are taken modulo $N_x$ and $N_y$, respectively. These operators correspond to discrete first-order derivatives in space and time, and the discrete Laplacian operator, respectively. We also set 
$$
[D_h \alpha]_{i,j}^{k}:= ((D_1 \alpha)^k_{i,j},(D_1 \alpha)^k_{i-1,j}, (D_2 \alpha)^k_{i,j},(D_2 \alpha)^k_{i,j-1})\in\R^4
$$
and, given $w:\mathbb{T}_{h}^{2}\to\R^4$, indexed as $(w^{(1)},w^{(2)},w^{(3)},w^{(4)})_{i,j}$, the discrete divergence operator is given by
\begin{equation*}
(\div_h w)_{i,j}:=(D_1w^{(1)})_{i-1,j}+(D_1w^{(2)})_{i,j}+(D_2w^{(3)})_{i,j-1}+(D_2w^{(4)})_{i,j}.
\end{equation*}
Let $\mathcal{K}=[0,+\infty)\times(-\infty,0]\times[0,+\infty)\times (-\infty,0]$ and, for every $z\in\R^{4}$, set $P_{\mathcal{K}}(z)=\left([z_1]_{+},[z_{2}]_{-},[z_{3}]_{+},[z_{4}]_{-}\right)$ for its orthogonal projection onto $\mathcal{K}$, where, for $r\in\R$, $[r]_{+}=\max\{r,0\}$ and $[r]_{-}=\min\{0,r\}$.

In~\cite{MR2679575}, the authors propose a finite difference scheme for $({\rm MFG}_{\mathbb{T}^{2}})$ which, in the case of the Hamiltonian~\eqref{power_hamiltonian}, takes the form:
\begin{equation}
\label{eq:mfg_torus_discrete}
\begin{aligned}
\left\{
\begin{aligned}
&-(D_t u)^k_{i,j} -\nu(\Delta_{h}u)_{i,j}^{k} +\frac{1}{q'}\left|P_{\mathcal{K}}(-[D_h u]_{i,j}^{k})\right|^{q'} =f(x_{i,j},m^{k+1}_{i,j})\,,\\
&(D_t m)^k_{i,j} -\nu(\Delta_{h}m)_{i,j}^{k+1}+\mbox{div}_{h}\left( \left|P_{\mathcal{K}}(-[D_h u]^{k})\right|^{q'-2}P_{\mathcal{K}}(-[D_h u]^{k})m^{k+1}\right)_{i,j}=0\,,\\
&m^{0}_{i,j}=m_{0}(x_{i,j}), \hspace{0.3cm} u^{M}_{i,j}=g(x_{i,j},m^{M}_{i,j})\quad\text{ in }\mathbb{T}^{2}_{h},
\end{aligned}
\right.
\end{aligned}
\tag{${\rm MFG}_{\mathbb{T}^{2}_{h}}$}
\end{equation}
where the first two equations hold in $\mathbb{T}^{2}_{h}\times T_{\dt}$. The convergence of solutions to~\eqref{eq:mfg_torus_discrete} towards solutions to $({\rm MFG}_{\mathbb{T}^{2}})$ has been established in~\cite{MR3097034} for classical solutions and in~\cite{MR3452251} for weak solutions.

For Neumann boundary conditions, we set the discretization parameters $h=(\dx,\dy)=((b-a)/(N_x+1),(d-c)/(N_y+1))$ such that the internal nodes and the boundary nodes are given by, respectively,
\begin{align*}
    \Omega_h&:= \{x_{i,j}=(a+i\dx,c+j\dy) \; ; \; i=1,2,\ldots,N_{x},\;j=1,2,\ldots,N_{y}\}, \\
    \partial\Omega_{h}&:= \{x_{0,j},x_{N_x+1,j},x_{i,0},x_{i,N_y+1} \; ; \; i=0,1,\ldots,N_{x}+1,\;j=0,1,\ldots,N_{y}+1\}.
\end{align*}
This definition ensures that the number of unknowns in the space--time grid remains the same irrespective of boundary conditions.
The finite difference operators in space are now defined only for the internal nodes, not modulo $N_x$ nor $N_y$.
The discretization of (MFG$_{\Omega})$ is equivalent to \eqref{eq:mfg_torus_discrete} with $\mathbb{T}^2_h$ replaced by $\Omega_h$ and the following additional constraints defined on $\partial\Omega_{h}$:
\begin{align*}
&-\nu(D_1 m)^{k+1}_{0,j}+\Phi(u)^k_{0,j,1}m^{k+1}_{0,j} +\Phi(u)^k_{1,j,2}m^{k+1}_{1,j}=0, \\
&\nu(D_1 m)^{k+1}_{N_x,j}-\Phi(u)^k_{N_x,j,1}m^{k+1}_{N_x,j} -\Phi(u)^k_{N_x+1,j,2}m^{k+1}_{N_x+1,j}=0, \\
&-\nu(D_2 m)^{k+1}_{i,0}+\Phi(u)^k_{i,0,3}m^{k+1}_{i,0} +\Phi(u)^k_{i,1,4}m^{k+1}_{i,1}=0, \\
&\nu(D_2 m)^{k+1}_{i,N_y}-\Phi(u)^k_{i,N_y,3}m^{k+1}_{i,N_y}-\Phi(u)^k_{1,N_y+1,4}m^{k+1}_{i,N_y+1}=0, \\
&(D_1u)_{0,j}^{k}=(D_1u)_{N_x,j}^{k}=(D_2u)_{i,0}^{k}=(D_2u)_{i,N_y}^{k}=0,
\end{align*}
where $\Phi(u)^k_{i,j,\cdot}$ denotes the entries of $\left( \left|P_{\mathcal{K}}(-[D_h u]^{k})\right|^{q'-2}P_{\mathcal{K}}(-[D_h u]^{k})\right)_{i,j}$ (out of 4).
We denote this discretized system by (MFG$_{\Omega_h})$.
The Neumann boundary conditions with mixed indices in the terms $\Phi(u) m$ are inspired by those implemented in \cite{achdou_mean_2020}. However, they use fictitious boundary nodes and set all (boundary) terms involving $\Phi(u)$ to zero using the boundary information about $u$. This nevertheless leads to the same system of equations up to relabelling of nodes.
Alternatively, avoiding a mismatch of indices in the terms $\Phi(u) m$, we obtain again that one entry of $\Phi(u)$ is out of bounds for each condition. There are two remedies for this problem: either introduce a variable $w=\Phi(u) m$ such that the underlying information becomes irrelevant, or introduce fictitious Neumann boundary conditions on $u$, as the fictitious boundary conditions set the four out-of-bounds terms to zero in that case.

\subsection{Variational formulation}\label{sec:Variational}
As noticed in~\cite{LasryLions07}, systems of the form~\eqref{MFGcont_torus}--(MFG$_{\Omega})$ admit a variational formulation. Define $b:\R \times \R^{d}\to \R \cup \{+\infty\}$ and $F$, $G:\mathbb{T}^{d}\times \R \to \R \cup \{+\infty\}$ as
\begin{align*}
\label{functions_definitions} 
b(m,w) &:= \left\{
\ba{ll}
\frac{1}{q}\frac{|w|^{q}}{m^{q-1}}
&\text{ if } m>0,\\
0 &\text{ if } (m,w)=(0,0), \\
+\infty &\text{ otherwise,}
\ea
\right.\\
F(x,m) &:= \left\{ \ba{ll} \int\limits_{0}^{m} f(x,m') \,dm' & \mbox{if } m\geq 0, \\
+\infty & \mbox{otherwise,}\ea \right. \\G(x,m) &:= \left\{ \ba{ll} \int\limits_{0}^{m} g(x,m') \,dm' & \mbox{if } m\geq 0, \\[4pt]
+\infty & \mbox{otherwise,}\ea \right.\,
\end{align*}
where $q$ is the conjugate exponent of $q'$, i.e., $q = q'/(q' -1)$.
System~\eqref{MFGcont_torus} is the optimality system of the variational problem:
\begin{equation}
\label{MFGvar}
\left\{
\begin{aligned}
&\inf\limits_{(m,w)} \int\limits_{0}^{T}\int\limits_{\mathbb{T}^{d}} b(m(x,t), w(x,t)) + F(x,m(x,t))  \ddx\ddt+ \int\limits_{\mathbb{T}^{d}} G(x,m(x,T)) \ddx& \\
&\mbox{s.t. } \hspace{0.3cm}\partial_t m -\nu \Delta m + \mbox{div}(w) =0\hspace{0.3cm}\mbox{in } \; \mathbb{T}^{d}\times [0,T], \\
&\hspace{0.9cm}m(\cdot,0) = m_0(\cdot)\hspace{0.2cm}\mbox{in } \mathbb{T}^{d},
\end{aligned}\right.\tag{VMFG$_{\mathbb{T}^{d}}$}
\end{equation}
and the system $({\rm MFG}_{\Omega})$ is the optimality system of (VMFG$_{\mathbb{T}^{d}}$) with ${\mathbb{T}^d}$ replaced by $\Omega$ and the following additional constraint
$$\left(\nu\nabla m-w\right)\cdot n=0 \;\text{ on } \;\partial\Omega\times [0,T].\\
$$
We denote this variational problem by (VMFG$_{\Omega})$.
Under the assumption that $f(x,\cdot)$ and $g(x,\cdot)$ are non-decreasing, the formulations \eqref{MFGvar}, (VMFG$_{\Omega})$ are shown to correspond to convex optimization problems and convex duality techniques can be applied to analyse \eqref{MFGcont_torus} and~(MFG$_{\Omega})$. See, e.g., \cite{Cardaliaguet_var,Benamou_Carlier_Santambrogio,Cesaroni_Cirant_var,Santambrogio_var} and the references therein for~\eqref{MFGcont_torus} or its ergodic version,  \cite{meszaros_silva_preprint_17} for a stationary version of~(MFG$_{\Omega})$, and~\cite{Gomes_Ricciardi} for the first-order case $\nu=0$ in~(MFG$_{\Omega})$. 

In the following, we consider the numerical approximations of the continuous optimization problems~\eqref{MFGvar} and~(VMFG$_{\Omega})$, the optimality conditions of which lead to \eqref{eq:mfg_torus_discrete} and (MFG$_{\Omega_h})$, respectively.
Define $\widehat{b}:\mathbb{R}\times\mathbb{R}^4 \to   \R\cup \{+\infty\}$ as
\begin{equation*}
\widehat{b}(m,w):=
\left\{
\ba{ll}
\frac{1}{q}\frac{|w|^{q}}{m^{q-1}},
&\text{ if } m>0, ~ w \in \mathcal{K},\\[6pt]
0, &    \text{ if } (m,w)=(0,0), \\[6pt]
+\infty, &\text{ otherwise},
\ea
\right.
\end{equation*}
and for $m:\mathbb{T}^2_h\times T_{\dt}\to\R$ and $w:\mathbb{T}^2_h\times (T_{\dt}\setminus\{T\})\to\R$ set
\begin{align*}
\B(m,w)&:=\sum_{i,j,k}\widehat b(m_{i,j}^k, w_{i,j}^k), \\
\F(m)&:=\sum_{i,j,k} F(x_{i,j},m_{i,j}^k) + \frac{1}{\dt}\sum_{i,j} G(x_{i,j},m_{i,j}^{N_t}), \\
\G(m,w)&:=\begin{bmatrix}
    m_{i,j}^0- m_0(x_{i,j}) \\
\ba{c}(D_tm)^0_{i,j} -\nu(\Delta_hm)^1_{i,j}  + (\div_h w)^{0}_{i,j} \\\vdots\\(D_tm)^{N_t-1}_{i,j} -\nu(\Delta_hm)^{N_t}_{i,j}  + (\div_h w)^{N_t-1}_{i,j}\ea
\end{bmatrix}.
\end{align*}
Then, problem~\eqref{eq:mfg_torus_discrete} is the optimality condition of the discrete variational problem (see, e.g., \cite{briceno-arias_implementation_2019}):
\begin{equation}
\label{MFGvarh}
\left\{
\begin{aligned}
\inf\limits_{(m,w)}\hspace{0.2cm}&\B(m,w)+\F(m)& \\
\mbox{s.t.} \hspace{0.3cm}&\G(m,w)=\bo,
\end{aligned}\right. \tag{VMFG$_{h}$}
\end{equation}
where $\bo$ stands for a zero vector of $(N_t + 1)N_xN_y$ entries. Applying the initial condition leaves $N_tN_xN_y$ unknowns for the mass and the same number for each component of $w$, leading to a total of $5N_tN_xN_y$ unknowns.

Similarly, problem~(MFG$_{\Omega_h})$ is the optimality condition of the variational problem \eqref{MFGvarh} with $\mathbb{T}^2_h$ replaced by $\Omega_h$ and the following additional constraints defined on $\partial\Omega_h$:
\begin{align*}
-\nu(D_1 m)^{k+1}_{0,j}+(w^{(1)})^k_{0,j} + (w^{(2)})^k_{1,j}&=0, \\
\nu(D_1 m)^{k+1}_{N_x,j}-(w^{(1)})^k_{N_x,j}-(w^{(2)})^k_{N_x+1,j}&=0, \\
-\nu(D_2 m)^{k+1}_{i,0}+(w^{(3)})^k_{i,0}+(w^{(4)})^k_{i,1}&=0, \\
\nu(D_2 m)^{k+1}_{i,N_y}-(w^{(3)})^k_{i,N_y}-(w^{(4)})^k_{i,N_y+1}&=0,
\end{align*}
applied at each time step, leading again to a total of $5N_tN_xN_y$ unknowns.
Note that the operator $\Delta_h$ is symmetric, hence if we introduce discrete Lagrange multipliers $u^{k}_{0,j},u^{k}_{N_x+1,j},u^{k}_{i,0},u^{k}_{i,N_y+1}$
associated with the Neumann boundary constraints and differentiate the discrete Lagrangian with respect to $m$, we obtain again the discrete Neumann boundary conditions for $u$.

The following sections are devoted to the efficient numerical solution to systems~\eqref{eq:mfg_torus_discrete} and (MFG$_{\Omega_h})$ through variational and  numerical linear algebra techniques. 

\section{A primal--dual method for the discretized variational mean field game} \label{sec:alg}
In this section we state a version of the primal--dual method introduced in \cite{MR2782122} to solve the finite-dimensional convex optimization problem \eqref{MFGvarh}. For a convex optimization problem of the form
\begin{equation*}
\min_{y\in\R^n}{\varphi(y)+\psi(y)},
\end{equation*}
	where $\varphi\colon\R^n\to \R \cup \{+\infty\}$ and $\psi\colon\R^n\to \R \cup 
\{+\infty\}$ are 
convex l.s.c. proper functions, the $(\ell +1)$-th Chambolle--Pock iteration reads
\begin{equation}
\label{pasosChamPock}
\begin{aligned}
x_{\ell+1}&:= 
\prox_{\sigma\psi^*}(x_{\ell}+ \sigma\tilde{y}_{\ell}),\\[6pt]
y_{\ell+1}&:=
\prox_{\tau\varphi}(y_{\ell}-\tau x_{\ell+1}),\\[6pt]
\tilde{y}_{\ell+1}&:= 
y_{\ell+1}+ \theta(y_{\ell+1}-y_{\ell}),
\end{aligned}\end{equation}
where  $\sigma>0$ and $\tau>0$ satisfy $\sigma\tau<1$, the parameter $\theta\in[0,1]$, and for a given a l.s.c. convex proper function 
$\phi:\R^n\to\, ]-\infty,+\infty]$, the proximal operator is defined as
\begin{align} \label{eq:def_prox_operator}
\prox_{\gamma\phi}x &:= 
\underset{y\in\R^n}{\mbox{argmin}}\left\{\phi(y)+\frac{|y-x|^2}{2\gamma}\right\}
\hspace{0.5cm} \forall \; x\in \R^n.
\end{align}
As in~\cite{briceno-arias_implementation_2019}, we approximate the solution to~\eqref{eq:mfg_torus_discrete} by considering the iterates~\eqref{pasosChamPock} applied to Problem~\eqref{MFGvarh} with  $y=(m,w)^\top$, $\varphi(m,w):=\B(m,w)+\F(m)$, and
\begin{align*}
\psi(m,w):=\begin{cases}
0 & \text{if}\;\G(m,w)=\bo,\\
+\infty & \text{otherwise}.
\end{cases}
\end{align*}
The first and second steps in~\eqref{pasosChamPock} require the computation of $\prox_{\sigma\psi^*}$ and $\prox_{\tau\varphi}$, respectively. Regarding the latter, it follows directly from~\eqref{eq:def_prox_operator} that 
\begin{equation*}
\left(\prox_{\tau\varphi}(m,w)\right)_{i,j}^{k}=\prox_{\tau\varphi_{i,j}}(m_{i,j}^{k},w_{i,j}^k), 
\end{equation*}
where $\varphi_{i,j}(m,w):=\widehat b(m,w)+ F(x_{i,j},m)$ for all $(m,w)\in\R\times\R^4$. Therefore, $\prox_{\tau\varphi}$ can be directly computed from the semi-explicit expressions, involving the solution of a scalar equation, for $\prox_{\tau\varphi_{i,j}}$ in~\cite[Corollary 3.1]{BAKS}. For the computation of $\prox_{\sigma\psi^*}$,
we begin by noting that the constraint $\G(m,w)=\bo$ can be conveniently expressed as the linear system
\begin{align}\label{ABC}
\bC(m,w)^\top:=[\bA|\bB](m,w)^\top=\bd\,,
\end{align}
where $\bA$ is an $(N_t+1)N_xN_y\times (N_t+1)N_xN_y$ matrix arising from the discretization of the term $(\partial_t-\nu\Delta)m$, $\bB$ is an $(N_t+1)N_xN_y\times 4(N_t+1)N_xN_y$ matrix associated to the discretization of $\mbox{\div}(w)$, and $\bd$ is the vector incorporating the initial datum $m_0(x)$.
Next, the identity (see, e.g., \cite[Section 24.2]{MR3616647}) 
\begin{align*}
\prox_{\sigma\psi*}=\Id-\sigma \prox_{\psi/\sigma}\circ(\Id/\sigma)=\Id-\sigma 
\prox_{\psi}\circ(\Id/\sigma),
\end{align*}
together with
 \begin{align*}
\prox_{\psi}\colon(m,w)\mapsto 
(m,w)^\top-\bC^\top(\bC\bC^\top)^{-1}(\bC(m,w)^\top-{\bd})
 \end{align*}
being determined as the solution of \eqref{ABC} with smallest (squared) $2$-norm distance to $(m,w)^\top$, leads to
\begin{align}\label{CC*}
\prox_{\sigma\psi*}\colon(m,w)\mapsto 
\bC^\top(\bC\bC^\top)^{-1}(\bC(m,w)^\top-\sigma\bd)\,.
\end{align}

We may now briefly state the structure of the matrices arising from the finite difference discretization of \eqref{ABC} at each Chambolle--Pock iteration. The matrices $\bA$ and $\bB$, arising from the terms $(\partial_t-\nu\Delta)m$ and $\mbox{\div}(w)$ respectively, are
\begin{equation*}
\bA = \left[\begin{array}{cccccc}
I & 0 & \cdots & \cdots & 0 & 0 \\
- \frac{1}{\dt} I & L & 0 & & & 0 \\
0 & - \frac{1}{\dt} I & L & \ddots & & \vdots \\
\vdots & \ddots & \ddots & \ddots & \ddots & \vdots \\
0 & & \ddots & - \frac{1}{\dt} I & L & 0 \\
0 & 0 & \cdots & 0 & - \frac{1}{\dt} I & L \\
\end{array}\right], \quad \bB = \left[\begin{array}{cccccc}
0 & 0 & \cdots & \cdots & 0 & 0 \\
B & 0 && & & 0 \\
0 & B & \ddots && & \vdots \\
\vdots & & \ddots & \ddots & & \vdots \\
0 & & & B & 0 & 0 \\
0 & 0 & \cdots & 0 & B & 0 \\
\end{array}\right].
\end{equation*}
Here, $I$ is the identity matrix, $B$ discretizes the divergence operator at each time step, and $L = \nu K + \frac{1}{\dt} I$ with $K$ discretizing the negative Laplacian operator. The matrices $L$ and $I$ are symmetric positive definite, with $K$ and $B B^\top$ symmetric positive semidefinite.

Upon forming $\bC\bC^\top$, which is equivalent to $\bA\bA^\top+\bB\bB^\top$ in the notation above, we wish to solve a system involving $\bC\bC^\top$, as required for \eqref{CC*}. We may eliminate the first ($N_xN_y \times N_xN_y$) block row, arising from the initial condition, and are then left with a reduced $N_tN_xN_y\times N_tN_xN_y$ system of the form $\bm{\mathcal A}\bx=\bb$, where
\begin{equation}\label{eqn:largeA}
\bm{\mathcal A} = \left[\begin{array}{cccccc}
C & - \frac{1}{\dt} L & 0 & \cdots & 0 & 0 \\
- \frac{1}{\dt} L & C + \frac{1}{(\dt)^2} I & - \frac{1}{\dt} L & \ddots & & 0 \\
0 & - \frac{1}{\dt} L & C + \frac{1}{(\dt)^2} I & - \frac{1}{\dt} L & & \vdots \\
\vdots & \ddots & \ddots & \ddots & \ddots & 0 \\
0 & & \ddots & - \frac{1}{\dt} L & C + \frac{1}{(\dt)^2} I & - \frac{1}{\dt} L \\
0 & 0 & \cdots & 0 & - \frac{1}{\dt} L & C + \frac{1}{(\dt)^2} I \\
\end{array}\right],
\end{equation}
where $C = L^2 + B B^\top$. As $\bm{\mathcal A}$ can be a huge-scale matrix, particularly upon fine discretizations of the space and time variables, it is particularly important to be able to accurately and efficiently solve linear systems of this form. In the next section, we therefore consider how to effectively precondition the
matrix $\bm{\mathcal A}$, with the goal that such a preconditioner may be applied within the \gls{cg} algorithm \cite{CG}.

\section{Parallel-in-time iterative solver} \label{sec:prec}
Depending on the algebraic properties of the matrix $\bm{\mathcal A}$ \eqref{eqn:largeA} arising at each \gls{cp} iteration, the inverse can be directly computed, or the system can be solved directly. However, both options are expensive computationally and storing the potentially dense decompositions required for an exact solve is costly. Instead, the resulting system can be solved using a preconditioned Krylov subspace method. Note that the matrix $\bm{\mathcal A}$ \eqref{eqn:largeA} is symmetric positive definite, hence we can use the \gls{cg} method \cite{CG}. However, as the system arising at each \gls{cp} iteration is of very large scale, and is expected to be ill-conditioned with a fine mesh, preconditioning, i.e., transforming the system into a convenient form for the chosen solver, becomes essential. A preconditioner used within \gls{cg} will also need to be symmetric positive definite.

A limited number of preconditioning strategies for the system \eqref{eqn:largeA} have been explored in the literature, most popularly multigrid methods and incomplete Cholesky factorization as in \cite{briceno-arias_implementation_2019}, which both show satisfactory performance, particularly for large $\nu$.
For systems similar to \eqref{eqn:largeA}, multigrid methods \cite{AchdouPerez2012} and a block-diagonalizable preconditioner using a temporal transformation and multigrid in space \cite{andreev_17} have been applied. However, the former method, although developed on the discrete level, does not involve the operator $BB^\top$, and the latter is a continuous operator preconditioning approach involving
a temporal discretization similar to continuous Galerkin time-stepping. We employ a \enquote{discretize-in-space-first} approach and deal with the exact matrices that arise.
Two further approaches have been explored during this study: (i) a cheap-to-apply preconditioner neglecting the matrix $BB^\top$ where Chebyshev semi-iteration is used to approximate the inverse action of $L$; (ii) a refined preconditioner specifically tailored to the regime of a small viscosity $\nu$, employing a matching strategy \cite{PW12,PW13} and an incomplete Cholesky decomposition of the positive semidefinite matrix $BB^\top$. Given its triangular structure, the terms within the decomposition of $BB^\top$ are then applied using Gauss--Seidel iteration within the matching strategy. Both preconditioners show satisfactory performance, with the second showing generally superior performance in numerical tests as it incorporates more information about the underlying system, and is more robust with respect to the parameter $\nu$. 
In this paper, we instead derive suitable preconditioners which can be applied in parallel for each time step by employing real-valued \glspl{dft} to block-diagonalize the huge-scale system in the time variable (often referred to as a ParaDiag approach), which we find to be a superior approach in practice.

The derivation of parallel-in-time methods for PDEs is an area of increasing interest within numerical mathematics; we refer to \cite{Gander2015,GWZ25} for surveys. Alongside approaches such as the parareal scheme \cite{LMT01,MT02}, other multiple shooting strategies \cite{BZ89,CP93,Kiehl94}, space--time multigrid approaches like Parallel Full Approximation Scheme in Space and Time (PFASST) \cite{EM12,BMS17} and multigrid reduction-in-time (MGRIT) \cite{FFKMS14,DFFKM21}, all of which may be analysed in a unified manner \cite{GLRS23}, an emerging and active area of research is to treat the resulting linear systems within the ParaDiag fashion \cite{MR08,GLWXZ21,KMZ23}. We employ this philosophy within the present work, as this enables the resolution of highly-challenging systems via fast trigonometric transforms and a block diagonal solver as described above.

Parallel-in-time preconditioning approaches frequently arise from recognizing underlying symmetries and repeating block structures, a well-established principle in preconditioning \cite{wathen_preconditioning_2015}.
To this end, note that the matrix $\bm{\mathcal A}$ has a uniform pattern of diagonal and off-diagonal matrices, apart from the first $N_xN_y\times N_xN_y$ block. Assuming suitable boundary, initial, and possibly final conditions, the matrix $\bm{\mathcal A}$ can be shown to be similar to a recursive linear algebra decomposition involving a combination of negative Laplacian operators. In the zero viscosity case, the system $\bm{\mathcal A}\bm{x}=\bb$ is similar to the system arising from a three-dimensional Poisson equation, which can be solved exactly using \glspl{dct}, as noticed by Papadakis et al.\ \cite{MR3158785}. Conveniently, the discrete negative Laplacian operator can be diagonalized by a decomposition of the form $V^{-1} \Lambda V$, where $V$ is a \gls{dft} matrix (or a linear combination of \glspl{dft} in the form of \glspl{dct} or \glspl{dst}) determined by suitable boundary and initial (and final) conditions, and $\Lambda$ is a diagonal matrix containing the eigenvalues of the negative Laplacian operator.
We extend the solution strategy presented in \cite{MR3158785} to the full \gls{mfg} with non-zero viscosity by constructing suitable positive definite preconditioners on the discrete level for periodic and Neumann boundary conditions. Nevertheless, we allow for other geometries and non-constant viscosity $\nu(x)$ when applying the parallel-in-time preconditioning technique, as it remains true that the matrices discretizing the spatial derivatives are repeating at every time step. We omit the final condition included in the optimal transport problem \cite{MR3158785} (which would otherwise lead to the \gls{mfg} planning problem \cite{MR2888257}) and treat any target state as a penalty term in the objective function.

The main advantage of this approach is the parallelization capability of the preconditioners. The system arising at each time step can be solved in parallel with a number of solution methods, including LU factorization, \gls{amg}, or by extending the recursive linear algebra technique to also include spatial decompositions.
The following sections describe two options for the decomposition in time, a special no-viscosity case, and lastly two possible spatial decompositions in the case of a rectangular grid.

\subsection{Block-diagonalizable preconditioner in time}
We start by devising a preconditioned approximation for the discretization in time. Recall that $L=\nu K + \frac{1}{\Delta t}I_{N_xN_y}$, where from now on the subscript of an identity matrix corresponds to its dimension. Hence, on the discrete level we have
\begin{align*}
    C=L^2 + BB^\top =\hat{C} + \frac{2\nu}{\Delta t} K + \frac{1}{(\Delta t)^2}I_{N_xN_y},
\end{align*}
where
\begin{align*}
    \hat{C}=\nu^2 K^2 + BB^\top.
\end{align*}
Letting $H=\frac{1}{\Delta t} (2\nu K + \frac{1}{\Delta t}I_{N_xN_y})$ and $\hat{L} = \frac{1}{\Delta t}L$, then $\bm{\mathcal A}$ in \eqref{eqn:largeA} is of the form
\begin{align}
    \bm{\mathcal A} = \begin{bmatrix}
        \hat{C} + H & - \hat{L} & 0 & \cdots & 0 & 0 \\
        - \hat{L} & \hat{C} + 2\hat{L} & - \hat{L} & \ddots & & 0 \\
        0 & -\hat{L} & \hat{C} + 2\hat{L} & -\hat{L} & \ddots & \vdots \\
        \vdots & \ddots & \ddots & \ddots & \ddots & 0 \\
        0 & & \ddots & - \hat{L} & \hat{C} + 2\hat{L} & -\hat{L} \\
        0 & 0 & \cdots & 0 & -\hat{L} & \hat{C} + 2\hat{L}
    \end{bmatrix}. \label{eqn:Adecom}
\end{align}
Note that we can approximate $H$ with a scalar multiple of $\hat{L}$ by either subtracting a term $\frac{\nu}{\Delta t}K$, or adding a term $\frac{1}{(\Delta t)^2} I_{N_xN_y}$, in the first diagonal block, corresponding to setting $l=1$ and $l=2$ in $H\approx l\hat{L}$, respectively. Equivalently, consider a preconditioner of the form
\begin{align}
    \bm{\mathcal P} = \begin{bmatrix}
        \hat{C} + l \hat{L} & - \hat{L} & 0 & \cdots & 0 & 0 \\
        - \hat{L} & \hat{C} + 2\hat{L} & - \hat{L} & \ddots & & 0 \\
        0 & -\hat{L} & \hat{C} + 2\hat{L} & -\hat{L} & \ddots & \vdots \\
        \vdots & \ddots & \ddots & \ddots & \ddots & 0 \\
        0 & & \ddots & - \hat{L} & \hat{C} + 2\hat{L} & -\hat{L} \\
        0 & 0 & \cdots & 0 & -\hat{L} & \hat{C} + 2\hat{L}
    \end{bmatrix}. \label{eqn:prec}
\end{align}
We present the following lemma for the spectral properties of the preconditioned matrix $\bm{\mathcal P}^{-1}\bm{\mathcal A}$, which provide a guarantee for the convergence properties of the preconditioned system in an ideal setting.

\begin{lemma} \label{thm:eig}
    Suppose $\bm{\mathcal A}$ and $\bm{\mathcal P}$ are defined as in \eqref{eqn:Adecom} and \eqref{eqn:prec}, respectively, with two spatial dimensions $N_x=N_y$. Then 1 is an eigenvalue of $\bm{\mathcal P}^{-1}\bm{\mathcal A}$ with algebraic multiplicity at least $(N_t-1)N_xN_y$, and the remaining $N_xN_y$ eigenvalues are of the form $1+\mu$, with $\mu$ being eigenvalues of the matrix
    $\bm{\mathcal P}^{-1}\bm{\mathcal H},$
    where
    \begin{align*}
    \bm{\mathcal H} = \begin{bmatrix}
        H - l\hat{L} & 0 & \cdots & 0 \\
        0 & 0 & \ddots & \vdots \\
        \vdots & \ddots & \ddots & 0 \\
        0 & \cdots & 0 & 0
    \end{bmatrix}.
\end{align*}
\end{lemma}
\begin{proof}
Consider the eigenvalue problem
\begin{align}
    \bm{\mathcal A}\bm{x}=\lambda \bm{\mathcal P}\bm{x}, \label{eqn:eigA}
\end{align}
where $\bm{x}$ is the eigenvector associated with eigenvalue $\lambda$. Note that $\bm{\mathcal A} = \bm{\mathcal P} + \bm{\mathcal H}$.
Hence the eigenvalue problem \eqref{eqn:eigA} becomes
\begin{align}
    \bm{\mathcal P}\bm{x} + \bm{\mathcal H}\bm{x}=\lambda \bm{\mathcal P}\bm{x}\phantom{\mu} \quad \iff \quad \phantom{\mu}(\lambda - 1)\bm{\mathcal P}\bm{x} &= \bm{\mathcal H} \bm{x} \nonumber \\
    \iff \quad \phantom{\bm{\mathcal P}(\lambda- 1)}\mu\bm{x} &= \bm{\mathcal P}^{-1}\bm{\mathcal H} \bm{x}, \label{eqn:subeig}
\end{align}
where $\mu=\lambda-1$.
Note that by the dimension of the problem that there are at most $N_tN_xN_y$ linearly independent non-trivial solutions $\bm{x}$ to the eigenvalue problem \eqref{eqn:subeig}. 
The matrix $\bm{\mathcal H}$ is made up of mostly zeros and has rank at most $N_xN_y$. Hence there are at least $(N_t-1)N_xN_y$ linearly independent solutions $\bm{x}$ which belong to the nullspace of the matrix in \eqref{eqn:subeig}. For these $(N_t-1)N_xN_y$ solutions, the eigenvalue problem \eqref{eqn:subeig} reduces to
$$(\lambda-1)\bm{x}={\bm 0},$$
which has the unique solution $\lambda=1$.
Therefore, the preconditioned matrix $\bm{\mathcal P}^{-1}\bm{\mathcal A}$ has at least $(N_t-1)N_xN_y$ eigenvalues equal to 1, and the remaining $N_xN_y$ are of the form $\lambda=1+\mu$, where $\mu$ is a solution to \eqref{eqn:subeig}.
\end{proof}

\begin{figure}[htb]
    \centering
    \begin{subfigure}[t]{.28\textwidth}
        \centering
        \includegraphics[width=\textwidth]{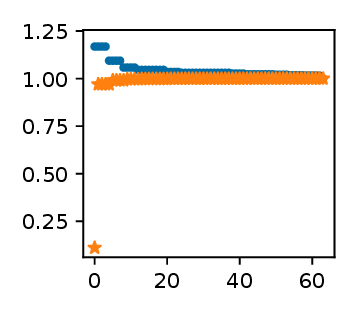}
        \caption{$\nu=1$}
    \end{subfigure}
    \begin{subfigure}[t]{.23\textwidth}
        \centering
        \includegraphics[width=\textwidth]{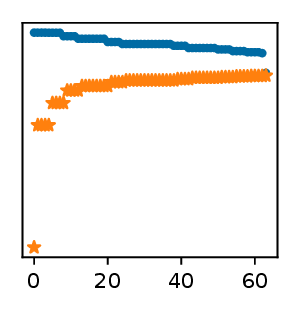}
        \caption{$\nu=0.1$}
    \end{subfigure}
    \begin{subfigure}[t]{.23\textwidth}
        \centering
        \includegraphics[width=\textwidth]{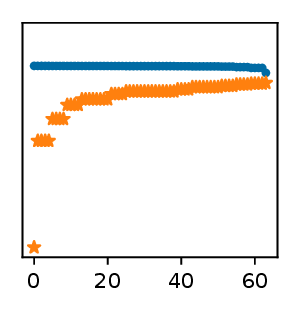}
        \caption{$\nu=0.01$}
    \end{subfigure}
    \begin{subfigure}[t]{.23\textwidth}
        \centering
        \includegraphics[width=\textwidth]{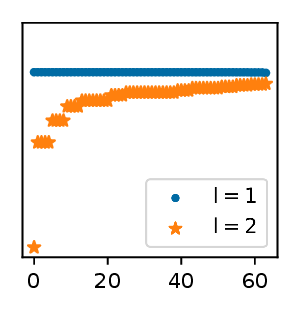}
        \caption{$\nu=0.001$}
    \end{subfigure}
    \caption{Non-unity eigenvalues of $\bm{\mathcal P}^{-1}\bm{\mathcal H}$ for different viscosity parameters $\nu$, and $N_x=N_y=N_t=8$. The horizontal axis contains the labels of the $N_xN_y$ eigenvalues, and the vertical axis displays the eigenvalues.} \label{fig:eigstime}
\end{figure}

Figure \ref{fig:eigstime} shows the $N_xN_y$ eigenvalues of $\bm{\mathcal P}^{-1}\bm{\mathcal H}$, i.e., those eigenvalues of $\bm{\mathcal P}^{-1}\bm{\mathcal A}$ not equal to 1, for the two preconditioning strategies and varying viscosity $\nu$. For $l=2$, one eigenvalue is relatively large in magnitude compared to the rest. However the eigenvalues remain bounded and consistent for all viscosity parameter values. For $l=1$, the spectrum of eigenvalues is small, with more non-unity eigenvalues in the particular case $\nu=0.1$. As the viscosity decreases, the eigenvalues converge to 1. This pattern can be shown to be consistent also in larger problem instances.

The preconditioner \eqref{eqn:prec} can be decomposed into the following sum of Kronecker products:
\begin{align*}
    \bm{\mathcal P} &= I_{N_t} \otimes \hat{C} + D_{tt} \otimes \hat{L},
\end{align*}
where $D_{tt}$ is a negative second-order finite difference operator
\begin{align*}
    D_{tt} = \begin{bmatrix}
        l & -1 & 0 & \cdots & 0 \\
        -1 & 2 & -1 &  \ddots & \vdots \\
        0 & -1 & 2 & \ddots & 0 \\
        \vdots & \ddots & \ddots & \ddots & -1 \\
        0 & \cdots & 0 & -1 & 2
    \end{bmatrix} \in \mathbb{R}^{{N_t} \times {N_t}}.
\end{align*}
Such an operator can be diagonalized by a decomposition of the type $D_{tt} = V^{-1} \Lambda_t V$, where $V$ is a matrix transform, either a type VIII \gls{dct} \cite{DCT} ($l=1$) or a type I \gls{dst} \cite{chan_sine_1997} ($l=2$), and $\Lambda_t$ is a diagonal matrix containing the eigenvalues of $D_{tt}$:
\begin{align}
    \lambda_k = -2\cos\left(\frac{\pi(k-1+\frac{l}{2})}{N_t+\frac{l}{2}}\right) +2, \quad k = 1, ..., N_t. \label{eqn:eigs}
\end{align}
The preconditioner can therefore be written as
\begin{align}
    \bm{\mathcal P} &= I_{N_t} \otimes \hat{C}  + (V^{-1} \Lambda_t V) \otimes \hat{L} \nonumber\\
    &= (V^{-1} \otimes I_{N_xN_y})(I_{N_t} \otimes \hat{C} + \Lambda_t \otimes \hat{L})(V \otimes I_{N_xN_y}). \label{eqn:precrec}
\end{align}
Note that changing from one preconditioner to the other ($l=1$, $l=2$) requires only that the chosen transformation $V$ corresponds to the correct eigenvalues $\lambda_k$ in \eqref{eqn:eigs}.

The inversion of the preconditioner \eqref{eqn:precrec} is described in Algorithm \ref{alg:recurt} in terms of its inverse action on an input vector $\bm{y}$.
We first reshape the input vector such that the matrix transformation $V$ can be applied to each column vector of dimension $N_t$, i.e., one application of $V$ is applied to the vector with precisely one corresponding entry from each of the subvectors $\underline{y}_k$, $k=1,...,N_t$.
The preconditioner \eqref{eqn:precrec} allows a parallel-in-time implementation when solving the system
\begin{align}
    (\hat{C} + \lambda_k \hat{L}) \underline{x}_k = \underline{y}_k, \label{eqn:paral}
\end{align}
for the inverse action of $\hat{C} + \lambda_k \hat{L}$ applied to an input vector $\underline{y}_k$ at each time step $k$. In the last step of the procedure in Algorithm \ref{alg:recurt}, the inverse transformation of $V$ is applied in the corresponding dimension, and the matrix solution reshaped into a column vector.

\begin{algorithm}[htb]
\caption{Application of the inverse action of $\bm{\mathcal P}$ in \eqref{eqn:precrec}} \label{alg:recurt}
\begin{algorithmic}[1]
\Require Input vector $\bm{y}=\begin{bmatrix}
\underline{y}^\top_1 & \underline{y}^\top_2 &\hdots &\underline{y}^\top_{N_t}\end{bmatrix}^\top\in \mathbb{R}^{N_tN_xN_y}$; preconditioner parameter $l\in\{1,2\}$
\Ensure Vector from inverse action of the preconditioner \eqref{eqn:precrec}, $\bm{x} = \bm{\mathcal P}^{-1} \bm{y}$
\State $\bm{x}\gets\begin{bmatrix}
\underline{x}^\top_1 & \underline{x}^\top_2 &\hdots &\underline{x}^\top_{N_t}\end{bmatrix}^\top\in \mathbb{R}^{N_tN_xN_y}$ \Comment{Define the output vector}
\State $\bm{y} \gets \begin{bmatrix}
    \underline{y}_1 & \underline{y}_2 & \hdots & \underline{y}_{N_t}
\end{bmatrix}$ \Comment{Reshape into an $N_xN_y \times N_t$ matrix}
\State $\bm{y} \gets V \bm{y}^\top$ \Comment{Transpose and apply column-wise matrix transformation}
\State $\bm{y} \gets \bm{y}^\top$ \Comment{Transpose}
\For{$k=1,...,N_t$}
\State $\underline{x}_k \gets (\hat{C} + \lambda_k \hat{L})^{-1}\underline{y}_k$ \Comment{Solve system \eqref{eqn:paral} for each column vector}
\EndFor
\State $\bm{x} \gets V^{-1} \bm{x}^\top$ \Comment{Transpose and apply inverse transformation}
\State $\bm{x} \gets \bm{x}^\top$ \Comment{Transpose}
\State $\bm{x} \gets \begin{bmatrix}
    \underline{x}_1^\top & \underline{x}_2^\top & \hdots & \underline{x}_{N_t}^\top
\end{bmatrix}^\top$ \Comment{Reshape back into a vector}
\end{algorithmic}
\Return $\bm{x}$ \Comment{Return the action applied to $\bm{y}$}
\end{algorithm}

When solving the system \eqref{eqn:paral}, there are a few options to consider in order to reduce runtime and memory. Firstly, the system can be solved exactly, using, e.g., an LU decomposition. This is a convenient method which is readily available with most commercial software. However, the runtime may be excessive for large problems. An incomplete Cholesky decomposition or an \gls{amg} method can be applied to reduce memory usage and runtime, but with the caveat that one obtains an approximate rather than exact solution.
In the case of a rectangular grid and constant viscosity, we can directly exploit the internal structure of the finite difference operators and consider a \gls{dft} or \gls{dct} decomposition in space. We describe this approach in Section \ref{subsec:recurs}.
Note also that in the no-viscosity case, the variational \gls{mfg} with $q'=2$ and $f\equiv0$ reduces to the Benamou--Brenier formulation of optimal transport without a final constraint \cite{benamou_computational_2000}, in which case the \gls{dct}-VIII decomposition in time, corresponding to the preconditioner \eqref{eqn:precrec} with $l=1$, is exact. We explain briefly in the following subsection this special case.

\subsection{The zero viscosity case} \label{subsec:novisc}
For small viscosity, the discrepancy between $H$ and $\hat{L}$ in the first diagonal block in $\bm{\mathcal A}$ \eqref{eqn:Adecom} vanishes, since
\begin{align*}
    \lim_{\nu \rightarrow 0} H = \frac{1}{(\Delta t)^2} I_{N_xN_y} = \lim_{\nu \rightarrow 0} \hat{L}.
\end{align*}
Therefore, $\bm{\mathcal P}$ with $l=1$ converges to $\bm{\mathcal A}$ as we decrease the viscosity, as subtracting the term involving $\nu K$ makes little to no difference. For small viscosities $0<\nu \ll 1$, we therefore expect this preconditioner to perform well. For $\nu=0$, an exact solution method can be employed. 
As observed in \cite{MR3158785}, the divergence operator $\text{div}$ applied to its adjoint operator, defined in terms of the $L^2$ inner product, is closely related to the negative Laplacian operator $-\Delta$, with some additional boundary terms. Also, the first order derivative $\partial_t$ with respect to time applied to its adjoint is closely related to the one-dimensional negative Laplacian operator $-\Delta$. Solving the system $\bm{\mathcal{\hat{A}}}\bm{z} = \bb$ without viscosity is therefore related to solving, on the continuous level, a three-dimensional Poisson equation: $$-\Delta z = b.$$
Let $D_{tt}$, $D_{xx}$, $D_{yy}$ denote the discrete negative Laplacian (or second derivative) operators. The discrete three-dimensional Laplacian operator $\bm{\mathcal{\hat{A}}}$ can be decomposed into the following sum of Kronecker products:
\begin{align*}
    \bm{\mathcal{\hat{A}}} &= \frac{1}{(\Delta t)^2}D_{tt} \otimes I_{N_x} \otimes I_{N_y} + \frac{1}{(\Delta x)^2}I_{N_t} \otimes D_{xx} \otimes I_{N_y}  + \frac{1}{(\Delta y)^2}I_{N_t} \otimes I_{N_x} \otimes D_{yy},
\end{align*}
allowing an exact solution using the \gls{dft} \cite{chan_circulant_1994} or a real-valued modification thereof in the form of \glspl{dst}/\glspl{dct}. In \cite{MR3158785}, both an initial and a final condition are assumed, leading to
\begin{align*}
    D_{tt} = \begin{bmatrix}
        1 & -1 & 0 & \cdots & 0 \\
        -1 & 2 & -1 &  \ddots & \vdots \\
        0 & -1 & \ddots & \ddots & 0 \\
        \vdots & \ddots & \ddots & 2 & -1 \\
        0 & \cdots & 0 & -1 & 1
    \end{bmatrix} \in \mathbb{R}^{{N_t} \times {N_t}},
\end{align*}
with Neumann boundary conditions in space. Including also a final condition, the problem can be easily solved with a three-dimensional \gls{dct}-II decomposition, which is a readily available solver in most common computing languages. We consider only an initial condition and include the option to apply a final penalty for any final misfit. Without a final condition, the lesser known \gls{dct}-VIII transformation can be used, however it is not available within popular scientific computing packages and requires a bespoke implementation.

When extending the idea of a recursive solution methods from optimal transport problems to variational \glspl{mfg}  with non-zero viscosity, no direct solution method exists as the system is no longer similar to a Poisson equation. Rather we get a fourth-order operator $KK^\top$ appearing in each diagonal block $\hat{C}$ which is notoriously ill-conditioned \cite{AchdouPerez2012}. Therefore, the system will need to be brought into a form that is easily solvable numerically.
We describe in the following section how we construct suitable decompositions in space for the discrete operators in $\bm{\mathcal A}$, arising from a \enquote{discretize-in-space-first} approach,
for periodic and Neumann boundary conditions on a rectangular grid and $\nu\geq 0$.

\subsection{Recursive preconditioner} \label{subsec:recurs}
We describe in this section the internal structure of the matrices $K$ and $B$ arising in $\hat{C}$ and $\hat{L}$ in \eqref{eqn:Adecom}, leading to a set of fully recursive preconditioners in both time and space.
Turning first to the divergence operator $B$, we treat the momentum variable $w$ separately in terms of its positive and negative components, as done in \cite{BAKS,briceno-arias_implementation_2019}.
The discrete divergence operator $B$ in one spatial dimension is of the form $\frac{1}{\Delta x}[D^1_x, D^2_x]$, where $D^1_x,D_x^2$ are (scaled) finite difference divergence operators. For periodic boundary conditions, we have $D^2_x=-(D^1_x)^\top$, where $D^1_x$ is defined as
$$
\begin{bmatrix}
        1 & 0 & \cdots & 0 & -1 \\
        -1 & 1 & 0 & & 0 \\
        0 & -1 & \ddots & \ddots & \vdots \\
        \vdots & \ddots & \ddots & 1 & 0 \\
        0 & \cdots & 0 & -1 & 1
    \end{bmatrix}.
$$
For the Neumann problem, the operators $D^1_x$, $D_x^2$ are (scaled) finite difference divergence operators with Dirichlet-type boundary conditions, respectively:
$$
\begin{bmatrix}
        1 & 0 & 0 & \cdots & 0 \\
        -1 & 1 & 0 &  & \vdots \\
        0 & -1 & \ddots & \ddots & 0 \\
        \vdots & \ddots & \ddots & 1 & 0 \\
        0 & \cdots & 0 & -1 & 0
    \end{bmatrix}, \quad \begin{bmatrix}
        0 & 1 & 0 & \cdots & 0 \\
        0 & -1 & 1 &  \ddots & \vdots \\
        0 & 0 & \ddots & \ddots & 0 \\
        \vdots & & \ddots & -1 & 1 \\
        0 & \cdots & 0 & 0 & -1
    \end{bmatrix}.
$$
Note that all the matrices are of size $N_x\times N_x$, while $\Delta x = (b-a)/N_x$ in the periodic case, and $\Delta x = (b-a)/(N_x+1)$ in the Neumann case.
Now, the discrete divergence operators multiplied with their transpose, $D^1_x (D^1_x)^\top$ and $D^2_x (D^2_x)^\top$, are both equal to, in the periodic and Neumann case, respectively,
$$
\begin{bmatrix}
    2 & -1 & 0 & \cdots & 0 & -1 \\
    -1 & 2 & -1 & \ddots & & 0 \\
    0 & -1 & 2 & \ddots & \ddots & \vdots \\
    \vdots & \ddots & \ddots & \ddots & -1 & 0 \\
    0 & & \ddots & -1 & 2 & -1 \\
    -1 & 0 & \cdots & 0 & -1 & 2
    \end{bmatrix}, \quad
    \begin{bmatrix}
    1 & -1 & 0 & \cdots & 0 & 0 \\
    -1 & 2 & -1 & \ddots & & 0 \\
    0 & -1 & 2 & \ddots & \ddots & \vdots \\
    \vdots & \ddots & \ddots & \ddots & -1 & 0 \\
    0 & & \ddots & -1 & 2 & -1 \\
    0 & 0 & \cdots & 0 & -1 & 1
    \end{bmatrix}.
$$
Adding these we have, after scaling, $BB^\top = \frac{1}{(\Delta x)^2}(D^1_x (D^1_x)^\top + D^2_x (D^2_x)^\top)$, or, respectively,
$$
\frac{2}{(\Delta x)^2}\begin{bmatrix}
    2 & -1 & 0 & \cdots & 0 & -1 \\
    -1 & 2 & -1 & \ddots & & 0 \\
    0 & -1 & 2 & \ddots & \ddots & \vdots \\
    \vdots & \ddots & \ddots & \ddots & -1 & 0 \\
    0 & & \ddots & -1 & 2 & -1 \\
    -1 & 0 & \cdots & 0 & -1 & 2
    \end{bmatrix}, \quad
\frac{2}{(\Delta x)^2}\begin{bmatrix}
    1 & -1 & 0 & \cdots & 0 & 0 \\
    -1 & 2 & -1 & \ddots & & 0 \\
    0 & -1 & 2 & \ddots & \ddots & \vdots \\
    \vdots & \ddots & \ddots & \ddots & -1 & 0 \\
    0 & & \ddots & -1 & 2 & -1 \\
    0 & 0 & \cdots & 0 & -1 & 1
    \end{bmatrix}.
$$
For this discretization, the discrete divergence operator multiplied with its transpose, $BB^\top$, is precisely two times the negative Laplacian operator with periodic or Neumann boundary conditions, i.e., $BB^\top=2K$. Hence we have
\begin{align}
    \hat{C}=\nu^2 K^2 + 2K. \label{eqn:precstead}
\end{align}

The discrete Laplacian operator with periodic boundary conditions is circulant and permits a \gls{dft} decomposition, and with Neumann boundary conditions the matrix permits a \gls{dct}-II decomposition.
Assume the discrete negative Laplacian operator $K$ permits one of the two decompositions, $K=\frac{1}{(\Delta x)^2}H^{-1}\Lambda_x H$, scaled by one spatial dimension, and with $\Lambda_x$ containing the eigenvalues of the negative Laplacian matrix. Since the Laplacian is symmetric, squaring this operator is equivalent to squaring only the diagonal matrix in the respective decomposition:
\begin{align*}
    K^2 &= \frac{1}{(\Delta x)^4}(H^{-1}\Lambda_x H)(H^{-1}\Lambda_x H)
    =\frac{1}{(\Delta x)^4}H^{-1}\Lambda_x^2 H.
\end{align*}
Therefore, the matrix $\hat{C}$ in \eqref{eqn:precstead} can be expressed in the convenient form
$$\hat{C}= \frac{\nu^2}{(\Delta x)^4} H^{-1}\Lambda_x^2 H + \frac{2}{(\Delta x)^2}H^{-1}\Lambda_x H=H^{-1}\left(\frac{\nu^2}{(\Delta x)^4}\Lambda_x^2+\frac{2}{(\Delta x)^2}\Lambda_x\right)H,$$
Note that this provides an easy extension from the optimal transport problem to the more general variational \glspl{mfg}, as we can set $\nu\geq0$ and retain the same decomposition.

Considering now two spatial dimensions, suppose $D_{xx}=H^{-1}\Lambda_xH$ and $D_{yy}=G^{-1}\Lambda_yG$, then the negative Laplacian operator is
\begin{align*}
    K&=\frac{1}{(\Delta x)^2}D_{xx}\otimes I_y + \frac{1}{(\Delta y)^2}I_x \otimes D_{yy} \\
    &= (I_{N_x}\otimes G^{-1})(H^{-1} \otimes I_{N_y})D (H \otimes I_{N_y})(I_{N_x}\otimes G),
\end{align*}
where
$$D = \frac{1}{(\Delta x)^2}\Lambda_x\otimes I_{N_y}  + \frac{1}{(\Delta y)^2}I_{N_x}\otimes\Lambda_y .$$
The matrix \eqref{eqn:precstead} can therefore be simply decomposed as follows:
\begin{align*}
    \hat{C} = (I_{N_x}\otimes G^{-1})(H^{-1} \otimes I_{N_y})\Lambda_1 (H \otimes I_{N_y})(I_{N_x}\otimes G),
\end{align*}
where
\begin{align}
    \Lambda_1 &= \nu^2 D^2 + 2D. \label{eqn:lambda1}
\end{align}
Similarly, we have
\begin{align*}
    \hat{L} &= \frac{1}{\Delta t}\left(\nu K + \frac{1}{\Delta t}I_{N_xN_y}\right) \\
    &= (I_{N_x}\otimes G^{-1})(H^{-1} \otimes I_{N_y})\Lambda_2(H \otimes I_{N_y})(I_{N_x}\otimes G),
\end{align*}
where
\begin{align}
    \Lambda_2 &= \frac{1}{\Delta t}\left(\nu D + \frac{1}{\Delta t}I_{N_x} \otimes I_{N_y}\right). \label{eqn:lambda2}
\end{align}
Hence the matrix \eqref{eqn:paral} arising at each time step $k$ is
\begin{align}
        \hat{C} + \lambda_k \hat{L} = (I_{N_x}\otimes G^{-1})(H^{-1} \otimes I_{N_y})(\Lambda_1 + \lambda_k \Lambda_2) (H \otimes I_{N_y})(I_{N_x}\otimes G).
        \label{eqn:paralrec}
\end{align}
This is a combination of matrix transforms and diagonal matrices, allowing an easy inversion.
The pseudocode for inverting the matrix \eqref{eqn:paralrec} is given in Algorithm \ref{alg:recurx}.
First, we define the eigenvalues of the negative Laplacian. Although the Laplacian with periodic or Neumann boundary conditions has a zero eigenvalue, a pertubation parameter is not needed in the time-dependent case due to the presence of $\frac{1}{(\Delta t)^2}I_{N_xN_y}$ in $\hat{L}$. Note also that the implementation allows mixed boundary conditions in the $x$ and $y$ directions, and extends easily to higher dimensions. The notation $\begin{bmatrix}
{\overline{y}}^\top_1 & \overline{y}^\top_2 &\hdots &\overline{y}^\top_{N_x}\end{bmatrix}^\top\in \mathbb{R}^{N_xN_y}$, refers to $N_x$ subvectors each of length $N_y$, i.e., $\overline{y}_i$ contains information about $\underline{y}$ at all spatial nodes in the second axis for index $i$ in the first axis. The diagonal matrix in \eqref{eqn:paralrec} can be reshaped into a convenient form $S$, leading to a computationally cheap entry-wise inversion on line 14. Line 10 refers to the addition of a constant to all entries of $S$. After reshaping $\underline{y}$ accordingly, a \gls{dft} or \gls{dct} is applied depending on the boundary conditions, along the $y$-direction first (line 12), and then the $x$-direction (line 13). After dividing by the eigenvalues, the inverse transformation is applied in order on lines 15--16. Lastly, the resulting matrix is reshaped into a vector.
\begin{algorithm}[h]
\caption{Application of the inverse action of $\hat{C}+\lambda_k\hat{L}$ in \eqref{eqn:paralrec}} \label{alg:recurx}
\begin{algorithmic}[1]
\Require Input vector $\underline{y}=\begin{bmatrix}
{\overline{y}}^\top_1 & \overline{y}^\top_2 &\hdots &\overline{y}^\top_{N_x}\end{bmatrix}^\top\in \mathbb{R}^{N_xN_y}$; boundary conditions; discretization parameters $\Delta x, \Delta y, \Delta t$; viscosity parameter $\nu$; eigenvalue $\lambda_k$ from the time decomposition
\Ensure Vector from inverse action of the matrix \eqref{eqn:paralrec}, $\underline{x} = (\hat{C}+\lambda_k\hat{L})^{-1} \underline{y}$
\State $\overline{v} \gets \begin{bmatrix}
0 & 1 &\hdots & N_x-1\end{bmatrix}^\top$ 
\If{periodic}
\State $\overline{\lambda}_x \gets -2\cos\left(\frac{2\pi \overline{v}}{N_x}\right) + 2$ \Comment{Define eigenvalues of the negative Laplacian}
\ElsIf{Neumann}
\State $\overline{\lambda}_x \gets -2\cos\left(\frac{\pi \overline{v}}{N_x}\right) +2$
\EndIf
\State Repeat lines 1--6 for $N_y$
\State $R \gets \frac{1}{(\Delta x)^2}\begin{bmatrix}
    \overline{\lambda}_x & \overline{\lambda}_x & \hdots & \overline{\lambda}_x 
\end{bmatrix}
+ \frac{1}{(\Delta y)^2}\begin{bmatrix}
    \overline{\lambda}_y & \overline{\lambda}_y & \hdots & \overline{\lambda}_y
\end{bmatrix}^\top
$ \Comment{Tile the eigenvalues}
\State $S \gets \nu^2R\odot R + 2R + \frac{\nu\lambda_k}{\Delta t}R$ \Comment{Obtain an $N_x\times N_y$ matrix}
\State $S \gets S + \frac{\lambda_k}{(\Delta t)^2}$ \Comment{Add constant to all entries}
\State $\underline{y} \gets \begin{bmatrix}
    \overline{y}_1 & \overline{y}_2 & \hdots & \overline{y}_{N_x}
\end{bmatrix}$ \Comment{Reshape into an $N_y \times N_x$ matrix}
\State $\underline{y} \gets G \underline{y}$ \Comment{Apply column-wise matrix transformation}
\State $\underline{y} \gets H \underline{y}^\top$ \Comment{Transpose and apply transformation}
\State $\underline{y} \gets \underline{y} \oslash S$ \Comment{Entry-wise division}
\State $\underline{y} \gets H^{-1} \underline{y}$ \Comment{Apply inverse transformation}
\State $\underline{y} \gets G^{-1} \underline{y}^\top$ \Comment{Transpose and apply inverse transformation}
\State $\underline{x} \gets \begin{bmatrix}
    \overline{y}_1^\top & \overline{y}_2^\top & \hdots & \overline{y}_{N_x}^\top
\end{bmatrix}^\top$ \Comment{Reshape back into a vector}
\end{algorithmic}
\Return $\underline{x}$ \Comment{Return the action applied to $\underline{y}$}
\end{algorithm}

Adding the decomposition in space to the time-dependent problem corresponds to adding layers to the recursive preconditioner \eqref{eqn:precrec}, that is
\begin{align}
    \bm{\mathcal P}
    = (V^{-1} \otimes H^{-1} \otimes G^{-1}) (I_{N_t} \otimes \Lambda_1 + \Lambda_t \otimes \Lambda_2)(V \otimes H\otimes G), \label{eqn:pstrec}
\end{align}
where $\Lambda_t$, $\Lambda_1$, $\Lambda_2$ are defined according to \eqref{eqn:eigs}, \eqref{eqn:lambda1}, \eqref{eqn:lambda2}, respectively.
Within an optimized parallel-in-time solver, we can replace line 6 in Algorithm \ref{alg:recurt} with the procedure Algorithm \ref{alg:recurx} and calculate the inverse action at each time step in parallel.

\section{Numerical results} \label{sec:num}
In this section, we showcase the performance of the set of preconditioners \eqref{eqn:prec} with different solution methods for the system \eqref{eqn:paral} arising at each time step.
We consider two (2+1)-dimensional test problems, one of which is the standard periodic problem \eqref{MFGcont_torus} with crowd aversion, and the other is the Neumann boundary problem (MFG$_{\Omega}$) with a terminal penalty on the density.
We generally consider grid resolutions in the range $N_x\in\{16,32,64,128\}$, and for $N_y=N_x$ and $N_t=8N_x$, the size of the system \eqref{eqn:largeA} is $8N_x^3$, hence increasing by a factor $2^3$ when the grid spacing is halved. The largest problem we consider has approximately 16.8 million unknowns arising at each \gls{cp} iteration when just accounting for the system \eqref{eqn:largeA}. We have in total five times this number of unknowns $(m,w)$ to be determined at each \gls{cp} iteration, i.e., up to 83.9m unknowns.

We describe here briefly the algorithm setup. 
In the accelerated \gls{cp} \cite{MR2782122}, the parameters are updated at each \gls{cp} iteration. In the first iteration, the primal and dual step sizes are both set to $\tau_0 = \sigma_0 = 1$.
At \gls{cp} iteration $\ell+1$, we first calculate the extrapolation parameter $$\theta_{\ell+1} = \frac{1}{\sqrt{1+2\gamma\tau_{\ell}}},$$ where the strong convexity constant $\gamma$ of the objective is set at the beginning of the algorithm, and hence update the primal and dual step sizes
$$\tau_{\ell+1} = \theta_{\ell+1}\tau_{\ell}, \quad
    \sigma_{\ell+1} = \frac{\sigma_{\ell}}{\theta_{\ell+1}}.$$
Note that the \gls{cp} parameters satisfy $\sigma > 0$, $\tau > 0$, $\theta \in [0, 1]$, and $\sigma \tau \leq 1$ at all \gls{cp} iterations. Although the accelerated \gls{cp} is only guaranteed to converge under strong convexity, we observed that a small constant $\gamma \sim 10^{-3}$ leads to good convergence in practice.
The \gls{cp} algorithm terminates when the following criterion is satisfied,
$$r = \|m_{\ell+1} - m_{\ell}\| \leq \texttt{CPtol},$$
i.e., when $r$, the change in $m$ at \gls{cp} iteration $\ell+1$ compared to the previous iteration, falls below a given tolerance \texttt{CPtol}. In general, we set $\texttt{CPtol}=10^{-4}\|m_0\|$, noting that $m_0$ is a probability distribution and hence cannot be identically zero. 
In test cases that use a preconditioner for the inner solve, we use \gls{cg} with tolerance set according to the rule $$\texttt{CGtol}=\min\{10^{-4}, \max\{10^{-6}, 10^{-4}r\}\}.$$

The solution method is implemented in \textit{Python 3.10} using \textit{SciPy} \cite{virtanen_scipy_2020} for standard linear algebra operations and solvers. The solution method, including the application of the preconditioner, is parallelized using the \textit{MPI4Py} \cite{mpi4py} library to ensure proper communication within the parallel infrastructure. We use 32 CPU cores allocated by an HPC cluster scheduler where different runs might use different processor models: \textit{Intel Xeon Gold 5218 (2.30 GHz)}, \textit{Intel Xeon Gold 6248 (2.50 GHz)}, or \textit{AMD EPYC 8534P (64-core)}.

\subsection{Problem 1: Crowd aversion with trigonometric cost terms} \label{subsec:period}
\begin{figure}[h]
    \centering
    \begin{subfigure}[b]{.23\textwidth}
        \centering
        \includegraphics[width=\textwidth]{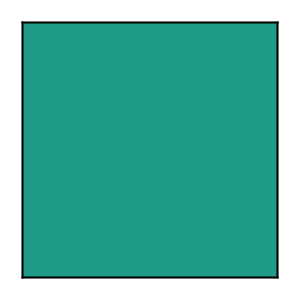}
        \caption{$t=0$}
    \end{subfigure}
    \begin{subfigure}[b]{.23\textwidth}
        \centering
        \includegraphics[width=\textwidth]{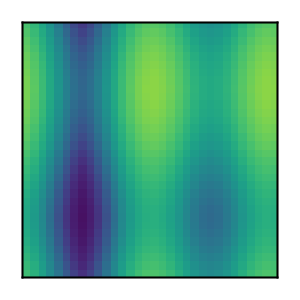}
        \caption{$t=0.1$}
    \end{subfigure}
    \begin{subfigure}[b]{.23\textwidth}
        \centering
        \includegraphics[width=\textwidth]{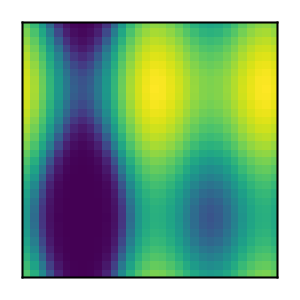}
        \caption{$t=0.5$}
    \end{subfigure}
    \begin{subfigure}[b]{.285\textwidth}
        \centering
        \includegraphics[trim=0 17 0 17,clip,width=\textwidth]{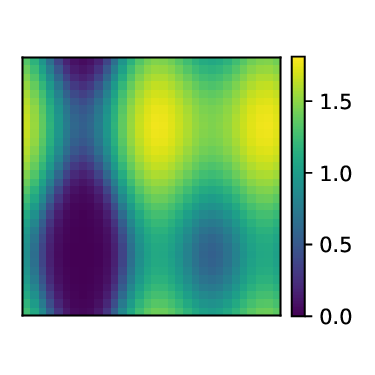}
        \caption{$t=1$}
    \end{subfigure}
    \caption{Solution to Problem 1 \eqref{MFGcont_torus}, with $N_x=N_y=32$, $N_t=5N_x$, $\nu=0.01$.} \label{fig:prob1}
\end{figure}

For this example, we simulate crowd aversion by setting
$$
f(x, y, m) = \frac12 \left(m^2-\sin(2\pi y)-\sin(2\pi x)-\cos(4\pi x)\right)
$$ 
and the Hamiltonian $H$ given by~\eqref{power_hamiltonian} with $q'=2$.
We set the initial density $m_0=1$, and $g\equiv 0$, and impose periodic boundary conditions and set $\gamma=5\cdot10^{-4}$.
The solution $m$ is shown at four different time steps in Figure \ref{fig:prob1}, validating the dynamic behaviour of the solution.
We observe similar behaviour to \cite{briceno-arias_implementation_2019} where the dynamic solution is close to the solution of the steady problem apart from the initial and final phase. Towards the end of the time interval, the cost function pushes the solution away from the stationary one, an observation known as the \textit{turnpike phenomenon} \cite{briceno-arias_implementation_2019}.

Table \ref{tbl:prob1all} shows the performance of the preconditioner \eqref{eqn:prec} for three problem sizes $N_tN_x^2$, with $N_t=8N_x$, and viscosity $\nu=0.01$. Columnwise, the results are given in terms of \gls{cp} iterations (\enquote{CP iter}), the average number of inner iterations (\enquote{Avg CG}) and CPU time in seconds (\enquote{Time [s]}) for two choices of preconditioner in time \eqref{eqn:prec} (\gls{dst}-I or \gls{dct}-VIII). The rows titled \enquote{LU} provide a benchmark in the form of an exact solve using LU decomposition with the \textit{SciPy SuperLU} solver pre-computed at the start of the \gls{cp} algorithm. Note that the \gls{dst}/\gls{dct} columns are irrelevant in this case. Evidently, the LU decomposition cannot be computed for the largest problem due to lack of memory with 32 GB (\enquote{NS}).
Next, we employ the parallel-in-time preconditioner \eqref{eqn:prec} and an LU solve for the system arising at each time step \eqref{eqn:paral} (\enquote{PinT + LU}).
We also employ a \textit{PyAMG Smoothed Aggregation} solver \cite{pyamgdoc, bell_pyamg_2022, vanek_algebraic_1996} for the term \eqref{eqn:paral}, with one (\enquote{PinT + 1 $V$-cycle}) and two (\enquote{PinT + 2 $V$-cycles}) $V$-cycles and no requirement for the relative residual norm tolerance.
Lastly, we employ a fully recursive preconditioner using a \gls{dft} decomposition in space \eqref{eqn:pstrec}. For all three approaches, the solvers for the systems arising at each time step \eqref{eqn:paral} are initialized at the start of the \gls{cp} algorithm and called at each \gls{cg} iteration.

The \gls{cp} iteration numbers are identical for all solution methods, and the \gls{cg} iteration numbers similarly remain stable for all problem sizes and solution methods, indicating robustness. The best performance in terms of inner iterations and CPU time are indicated in bold for each problem size. Clearly, the parallel-in-time preconditioner with a \gls{dct} decomposition in time and \gls{dft} in space offers the best results, as noted in bold. Note that the inner LU decomposition and the recursive decomposition are exact (at least in exact arithmetic), hence we obtain the same inner iteration numbers. Using a multigrid method for the system arising at each timestep, the inner iteration numbers increase slightly, however remaining contained with at most 36 inner iterations on average for the \gls{dst} decomposition in time, and 13 for the \gls{dct} decomposition. Using the recursive linear algebra method nevertheless leads to a decrease in both inner iteration numbers and CPU time. The problem is solved for 16.8m variables using the recursive preconditioner, which remains robust with respect to inner iteration numbers and CPU time.

\begin{table}[h]
    \centering
    \caption{Numerical results for the parallel-in-time preconditioner with various types of inner solves vs.\ an exact outer solve, for a periodic problem with viscosity $\nu=0.01$.} \label{tbl:prob1all}
    \small
    \begin{tabular}{l r r r r r r r}
    \toprule
    & $N_x$ & Prob.\ size & CP iter & Avg CG & Time [s] & Avg CG & Time [s] \\
    \cmidrule(lr){5-6} \cmidrule(lr){7-8}
    & & & & \multicolumn{2}{c}{(DST)} & \multicolumn{2}{c}{(DCT)} \\
    \midrule
    LU           & 16  & 32\,768     & 77   & --   & 9.79      & --   & -- \\
                 & 32  & 262\,144    & 67   & --   & 611.38    & --   & -- \\
                 & 64  & 2\,097\,152 & NS   & --   & NS        & --   & -- \\
    \midrule
    PinT          & 16  & 32\,768     & 77   & 11.8 & 10.78     & 2.4  & 29.14 \\
    + LU          & 32  & 262\,144    & 67   & 17.4 & 135.18    & 2.4  & 63.35 \\
                  & 64  & 2\,097\,152 & 70   & 26.1 & 325.17    & 2.4  & 266.22 \\
    \midrule
    PinT          & 16  & 32\,768     & 77   & 19.2 & 12.90      & 6.5  & 8.20 \\
    + 1 $V$-cycle & 32  & 262\,144    & 67   & 26.0 & 67.08     & 8.8  & 32.40 \\
                  & 64  & 2\,097\,152 & 70   & 35.7 & 326.67    & 12.9 & 283.13 \\
    \midrule
    PinT          & 16  & 32\,768     & 77   & 17.6 & 13.69     & 4.3  & 8.13 \\
    + 2 $V$-cycles& 32  & 262\,144    & 67   & 24.6 & 54.71     & 6.1  & 35.24 \\
                  & 64  & 2\,097\,152 & 70   & 34.1 & 367.92    & 8.9  & 268.77 \\
    \midrule
    PinT          & 16  & 32\,768     & 77   & 11.8 & 10.21     & 2.4  & 6.27 \\
    + recursive   & 32  & 262\,144    & 67   & 17.8 & 38.22     & 2.4  & 31.51 \\
                  & 64  & 2\,097\,152 & 70   & 26.1 & 268.66    & 2.4  & 234.78 \\
                  & 128 & 16\,777\,216& 68   & 38.5 & 2\,192.06   & 2.4  & 1832.57 \\

    \bottomrule
\end{tabular}
\end{table}

The performance of the fully recursive preconditioner is shown for different viscosities in Table \ref{tbl:prob1visc} in terms of average performance (\enquote{Avg DST/DCT perf}), i.e.\ the average \gls{cg} iteration numbers across all \gls{cp} iterations (\enquote{CG iter}) and the average CPU time for one \gls{cg} iteration (\enquote{Time [s]}). Both preconditioners show robustness with respect to the viscosity and problem sizes, and the CPU times are again similar for both. While there is an increase in outer iterations for smaller viscosity as the problem becomes harder to solve,
the \gls{cg} iteration numbers decrease. This is particularly the case for the \gls{dct} decomposition, as noted in bold, as it becomes exact for zero viscosity. Nevertheless, the \gls{dst} preconditioner dominates in terms of CPU time per \gls{cg} iteration, also noted in bold. This is partly due to the custom \gls{dct}-VIII implementation, which introduces an additional computational overhead. Performance could likely be improved with a more optimized version, such as a \textit{C++} implementation. As shown in Table~\ref{tbl:prob1all}, this highlights a trade-off: although \gls{dct} may reduce iteration counts, the readily available \gls{dst} preconditioner offers comparable total runtime due to its more efficient execution and smaller memory requirement.

\begin{table}[h]
    \centering
    \caption{Numerical results for the recursive parallel-in-time preconditioner, for a periodic problem with different viscosities.} \label{tbl:prob1visc}
    \small
    \begin{tabular}{l r r r r r r r}
    \toprule
    & & & & \multicolumn{2}{c}{Avg DST perf} & \multicolumn{2}{c}{Avg DCT perf} \\
    \cmidrule(lr){5-6} \cmidrule(lr){7-8}
    $\nu$ & $N_x$ & Prob.\ size & CP iter & CG iter & Time [s] & CG iter & Time [s] \\
    \midrule
    1            & 16 & 32\,768     & 20 & 7.4  & \textbf{0.01} & \textbf{4.7} & 0.03 \\
                 & 32 & 262\,144    & 20 & 8.9  & \textbf{0.06} & \textbf{4.9} & 0.11 \\
                 & 64 & 2\,097\,152 & 20 & 11.6 & \textbf{0.31} & \textbf{4.8} & 1.37 \\
    \midrule
    0.1          & 16 & 32\,768     & 21 & 13.9 & \textbf{0.02} & \textbf{3.9} & 0.04 \\
                 & 32 & 262\,144    & 21 & 19.5 & \textbf{0.09} & \textbf{3.9} & 0.13 \\
                 & 64 & 2\,097\,152 & 21 & 27.4 & \textbf{0.14} & \textbf{4.0} & 1.60 \\
    \midrule
    0.01         & 16 & 32\,768     & 77 & 11.8 & \textbf{0.01} & \textbf{2.4} & 0.03 \\
                 & 32 & 262\,144    & 67 & 17.8 & \textbf{0.03} & \textbf{2.4} & 0.20 \\
                 & 64 & 2\,097\,152 & 70 & 26.1 & \textbf{0.15} & \textbf{2.4} & 1.40 \\
    \midrule
    0.001        & 16 & 32\,768     & 69 & 13.1 & \textbf{0.03} & \textbf{2.0} & 0.15 \\
                 & 32 & 262\,144    & 56 & 17.8 & \textbf{0.03} & \textbf{2.0} & 0.24 \\
                 & 64 & 2\,097\,152 & 72 & 23.3 & \textbf{0.16} & \textbf{2.0} & 1.69 \\
    \bottomrule
\end{tabular}
\end{table}

\subsection{Problem 2: Gaussian final target state}
\begin{figure}[htb]
    \centering
    \begin{subfigure}[b]{.23\textwidth}
        \centering
        \includegraphics[width=\textwidth]{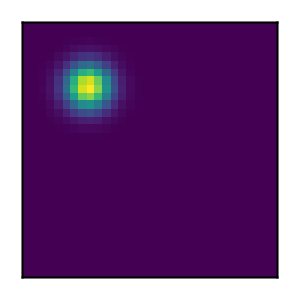}
        \caption{$t=0$}
    \end{subfigure}
    \begin{subfigure}[b]{.23\textwidth}
        \centering
        \includegraphics[width=\textwidth]{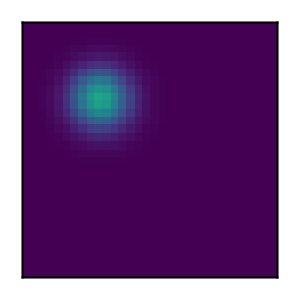}
        \caption{$t=0.1$}
    \end{subfigure}
    \begin{subfigure}[b]{.23\textwidth}
        \centering
        \includegraphics[width=\textwidth]{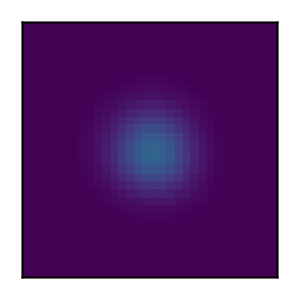}
        \caption{$t=0.5$}
    \end{subfigure}
    \begin{subfigure}[b]{.285\textwidth}
        \centering
        \includegraphics[trim=0 17 0 17,clip,width=\textwidth]{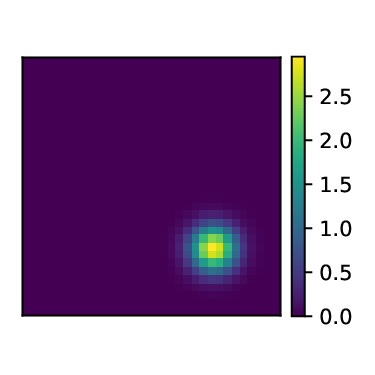}
        \caption{$t=1$}
    \end{subfigure}
    \caption{Solution to Problem 2 (MFG$_{\Omega}$) with a Gaussian target density, and $N_x=N_y=32$, $N_t=5N_x$, $\nu=0.01$, and $\beta\Delta t=10^{-5}$.} \label{fig:prob2}
\end{figure}

The second example we implement is the Neumann problem (MFG$_{\Omega}$) with $\Omega=[-1/2,1/2]\times[-1/2,1/2]$. We set a final penalty $g(x,y,m)=\frac{1}{\beta}(m-\bar{m})$, where $\bar{m}$ is the target density, and $\beta$ is a parameter determining the importance of reaching the target state, with small $\beta$ promoting convergence to the target. We set $\beta=10^{-3}$, noting that the term $g$ on the discrete level is weighted by $1/(\beta\Delta t)$, and the initial and target densities are Gaussian functions set to, respectively,
$$m_0=3\exp^{-2^7\left(\left(x+\frac{1}{4}\right)^2+\left(y-\frac{1}{4}\right)^2\right)}, \quad \bar{m}=3\exp^{-2^7\left(\left(x-\frac{1}{4}\right)^2+\left(y+\frac{1}{4}\right)^2\right)}.$$
We also consider the quadratic Hamiltonian~\eqref{power_hamiltonian} with $q'=2$, and set $f(x,y,m) \equiv 0$.
The \gls{cp} parameter $\gamma$ is set to $\gamma=5\cdot10^{-3}$.
The solution is shown in Figure \ref{fig:prob2} for viscosity $\nu=0.01$.

We repeat the same test as in Table \ref{tbl:prob1visc}, now for the Neumann problem, using an exact \gls{dct}-II decomposition in space \eqref{eqn:pstrec}.
This problem is significantly harder to solve, hence more outer iterations are needed to solve the problem and longer CPU times. Note, however, that the number of inner iterations is very similar to those shown for the periodic problem in Table \ref{tbl:prob1visc}, indicating robustness of the recursive preconditioners. The largest system we apply the preconditioner to has 2.1m variables, with 5 times this many unknowns at each \gls{cp} iteration.

\begin{table}[h]
    \centering
    \caption{Numerical results for the fully recursive parallel-in-time preconditioner, applied to the Neumann problem with different viscosities.} \label{tbl:prob2visc}
    \small
    \begin{tabular}{l r r r r r r r}
    \toprule
    $\nu$ & $N_x$ & Prob.\ size & CP iter & Avg CG & Time [s] & Avg CG & Time [s] \\
    \cmidrule(lr){5-6} \cmidrule(lr){7-8}
    & & & & \multicolumn{2}{c}{(DST)} & \multicolumn{2}{c}{(DCT)} \\
    \midrule
1 & 16 & 32\,768 & 109 & 11.5 & 9.59 & 5.2 & 6.89 \\
 & 32 & 262\,144 & 118 & 15.3 & 46.49 & 5.4 & 28.49 \\
 & 64 & 2\,097\,152 & 124 & 20.6 & 252.98 & 5.0 & 198.45 \\
 \midrule
0.1 & 16 & 32\,768 & 476 & 22.0 & 54.98 & 4.2 & 23.90 \\
 & 32 & 262\,144 & 549 & 30.3 & 233.01 & 4.0 & 139.83 \\
 & 64 & 2\,097\,152 & 674 & 40.5 & 1\,789.48 & 4.0 & 1\,219.40 \\
 \midrule
0.01 & 16 & 32\,768 & 354 & 21.6 & 38.84 & 3.8 & 16.33 \\
 & 32 & 262\,144 & 410 & 29.2 & 154.15 & 3.8 & 83.00 \\
 & 64 & 2\,097\,152 & 532 & 40.3 & 1\,058.02 & 3.6 & 641.61 \\
 \midrule
0.001 & 16 & 32\,768 & 351 & 20.3 & 33.68 & 2.6 & 16.08 \\
 & 32 & 262\,144 & 404 & 25.0 & 132.27 & 2.4 & 61.53 \\
 & 64 & 2\,097\,152 & 592 & 30.8 & 949.92 & 2.3 & 670.74 \\
    \bottomrule
    \end{tabular}
\end{table}

Next, we repeat the same test using an \gls{amg} solver for the system arising at each time step, shown in Table \ref{tbl:prob2viscamg}.
The \gls{amg} method performs very well for small viscosity but relatively poorly for large viscosity, whereas the recursive preconditioner shows stable performance for any viscosity in Table \ref{tbl:prob2visc}. Nevertheless, the spatial recursive linear algebra only works for a rectangular domain and a finite differences discretization, whereas \gls{amg} generalizes to any irregular domain.

\begin{table}[h]
    \centering
    \caption{Numerical results for the parallel-in-time preconditioner and an \gls{amg} inner solver, applied to the Neumann problem with different viscosities.} \label{tbl:prob2viscamg}
    \small
    \begin{tabular}{l r r r r r r r}
    \toprule
    $\nu$ & $N_x$ & Prob.\ size & CP iter & Avg CG & Time [s]  & Avg CG & Time [s] \\
    \cmidrule(lr){5-6} \cmidrule(lr){7-8}
    & & & & \multicolumn{2}{c}{(DST)} & \multicolumn{2}{c}{(DCT)} \\
    \midrule
    1     & 16      & 32\,768     & 109 & 50.0    & 38.38   & 48.6    & 42.87   \\
          & 32      & 262\,144    & 118 & 95.7    & 184.93  & 99.9    & 197.13  \\
          & 64      & 2\,097\,152 & 122 & 215.3   & 1\,569.09 & 220.5 & 1\,877.07 \\
    \midrule
    0.1   & 16      & 32\,768     & 476 & 27.7    & 114.81  & 16.3    & 73.47   \\
          & 32      & 262\,144    & 549 & 41.2    & 495.64  & 33.3    & 396.83  \\
          & 64      & 2\,097\,152 & 674 & 75.9    & 3\,796.97 & 63.9  & 3\,984.10 \\
    \midrule
    0.01  & 16      & 32\,768     & 354 & 28.8    & 76.75   & 6.0     & 25.41   \\
          & 32      & 262\,144    & 410 & 38.1    & 291.08  & 7.5     & 106.54  \\
          & 64      & 2\,097\,152 & 532 & 46.6    & 1\,942.38 & 9.4   & 961.39  \\
    \midrule
    0.001 & 16      & 32\,768     & 351 & 28.4    & 76.67   & 4.7     & 20.90    \\
          & 32      & 262\,144    & 404 & 36.0    & 276.09  & 5.1     & 82.59   \\
          & 64      & 2\,097\,152 & 592 & 44.3    & 1\,929.73 & 5.3   & 751.08  \\
    \bottomrule
\end{tabular}
\end{table}

\subsection{Parallelization}
To test the performance of the solution method with respect to scalability and parallelization, we perform a strong and weak scaling test. Firstly, we run three fixed size problems to see if increasing the number of cores decreases the runtime. Secondly, we fix the proportional problem size in each core and increase the number of cores to see if the runtime remains constant. All parameters are set to the same as mentioned previously.

The parallelization can be described as follows: Generally all vectors are divided into temporal subvectors. The number of subvectors that need to be defined for each core is the total number of time steps divided by the number of cores, currently implemented without a possible remainder, i.e.\ $N_t$ must be divisible by the number of cores. Various subvectors need to be communicated across cores due to the time-stepping in \eqref{eqn:largeA} involving the previous and next time steps. Similarly, there are special conditions for the subvectors involving the first and final time step. Matrix and vector multiplications are only performed for those matrices and vectors specified for each core. The \gls{cg} is parallelized by computing matrix and vector multiplications and vector norms separately and gathering them whenever necessary. For the preconditioner, the \gls{dst} or \gls{dct} transformation is applied in parallel by dividing the input vector of Algorithm \ref{alg:recurt} into \textit{spatial} subvectors and hence performing the transformation.
The vector is then redistributed in the time dimension, and the solution of \textit{temporal} subproblems \eqref{eqn:paral} is performed parallel-in-time. Lastly, the inverse transformation is applied spatially in parallel, and the output from \gls{cg} is redistributed evenly to the cores in order to apply the remaining steps of \gls{cp} (involving the proximal operators and the error norm) parallel-in-time. For the following scaling tests, we use a \gls{dct} decomposition in time and run exactly 100 iterations of the \gls{cp} algorithm.

\begin{figure}[h]
    \centering
    \includegraphics[width=0.6\linewidth]{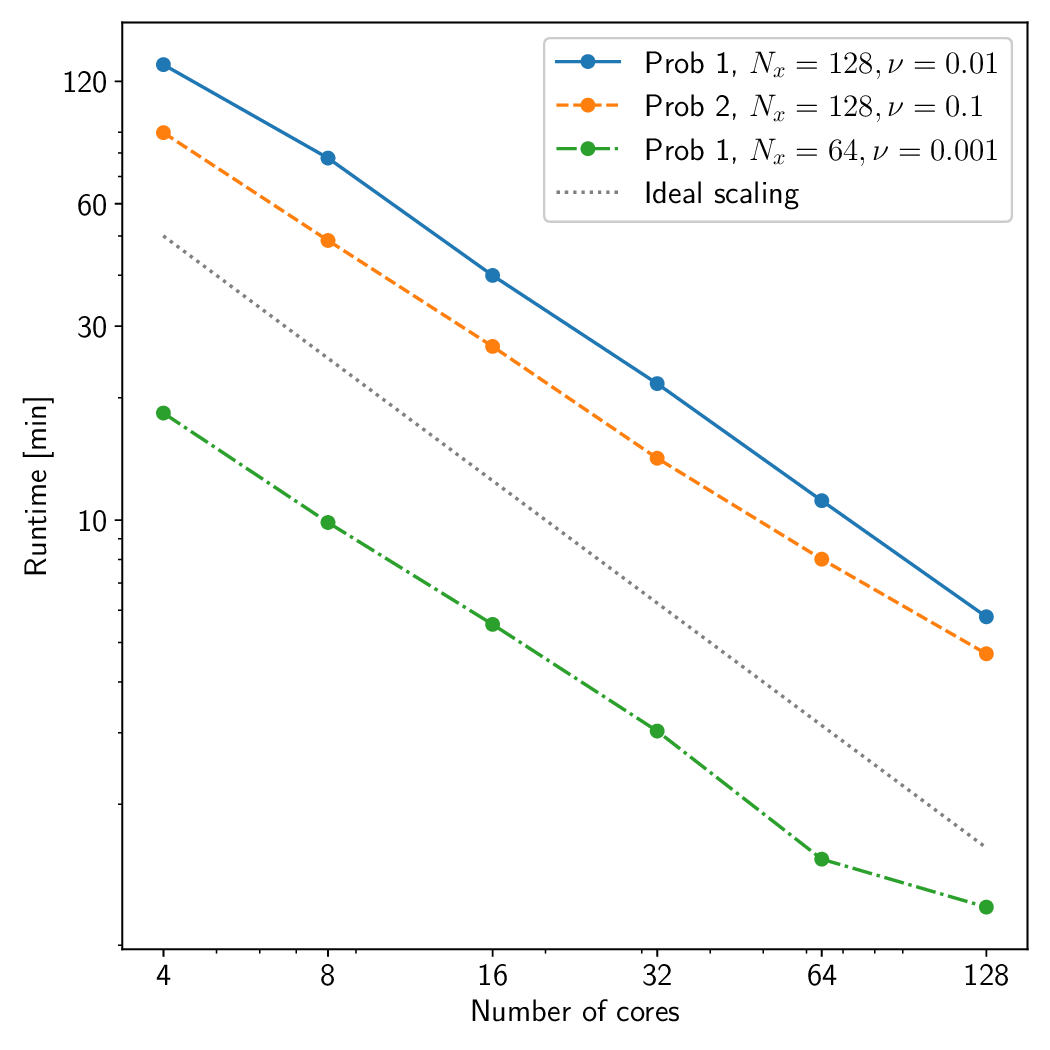}
    \caption{Strong scaling of the preconditioned method on a log--log scale, for three problem configurations.}
    \label{fig:paralstrong}
\end{figure}

For the strong scaling test, we set $N_y=N_x$ and $N_t=4N_x$ and run three test problem configurations: (i) Problem 1, $N_x=128$, $\nu = 0.01$; (ii) Problem 2, $N_x=128$, $\nu = 0.1$; and (iii) Problem 1, $N_x=64$, $\nu = 0.001$.
The two larger problems are allocated 8 GB of memory per core and the smaller one 6 GB per core, although actual usage showed that only 3 GB/core (1 GB/core for the smaller problem) was sufficient when using 4 cores, due to reduced communication buffer requirements.
Figure \ref{fig:paralstrong} shows the strong scaling performance of the recursive preconditioner.
Evidently, the runtime is approximately halved whenever double the amount of cores are available to the solver.
The slight runtime overhead can be explained by increased communication and synchronization costs between processors as the problem size per processor decreases.

Table \ref{tbl:weak} shows the weak scaling of the recursive \gls{dct} preconditioner when the size of the subproblem sent to each core remains constant. The spatial refinement is set to either $N_x=N_y=32$ or $N_x=N_y=64$, and the viscosity is set to $\nu=0.01$. The number of time steps $N_t$ is then set to 8 or 16 times the number of cores used for the smaller or larger problem, respectively. The two configurations are repeated for both spatial decompositions, i.e., the periodic (\enquote{Prob 1}) and the Neumann problem (\enquote{Prob 2}). Evidently, the runtime remains fairly constant across all problem instances. For the coarser grid, the runtime doubles in the worst case when the problem size increases by a factor 32. For the finer grid, the runtime remains relatively unchanged when the problem size increases by a factor 32, indicating solid scalability, but the runtime is inevitably subject to randomness due to fluctuating load and resource contention on the shared supercomputing cluster. Note additionally that the runtime overhead of the Neumann problem shown in Table \ref{tbl:prob2visc} is due to the number of \gls{cp} iterations, and not the preconditioning, as Table \ref{tbl:weak} shows shorter runtime per \gls{cp} iteration for the Neumann problem.

\begin{table}[h]
    \centering
    \caption{Weak scaling of the preconditioned method with 100 \gls{cp} iterations and $\nu=0.01$, for varying $N_t$ and two spatial configurations.}
    \label{tbl:weak}
    \small
    \begin{tabular}{lrrrrr}
        \toprule
        \# Cores & $N_t$ & Prob.\ size & \multicolumn{2}{c}{Runtime [s]} \\
        \cmidrule(lr){4-5}
        & & & Prob 1 & Prob 2 \\
        \midrule
        \multicolumn{5}{c}{\textbf{$N_x=N_y=32$}} \\
        \midrule
        4   & 32   & 32\,768     & 39.65   & 16.21  \\
        8   & 64   & 65\,536     & 40.28   & 22.87  \\
        16  & 128  & 131\,072    & 44.15   & 24.85  \\
        32  & 256  & 262\,144    & 46.01   & 26.03  \\
        64  & 512  & 524\,288    & 45.08   & 31.36  \\
        128  & 1024  & 1\,048\,576    & 50.59   & 34.42  \\
        \midrule
        \multicolumn{5}{c}{\textbf{$N_x=N_y=64$}} \\
        \midrule
        4   & 64   & 262\,144     & 414.07  & 144.60 \\
        8   & 128  & 524\,288     & 447.53  & 135.67 \\
        16  & 256  & 1\,048\,576  & 414.32  & 141.18 \\
        32  & 512  & 2\,097\,152  & 428.67  & 142.93 \\
        64  & 1024 & 4\,194\,304  & 408.46  & 163.20 \\
        128  & 2048 & 8\,388\,608 & 422.73  & 162.80 \\
        \bottomrule
    \end{tabular}
\end{table}

\section{Concluding remarks}
In this work, we introduced a novel preconditioning strategy for the efficient solution of linear systems arising in primal--dual methods applied to finite difference discretizations of time-dependent variational \gls{mfg} systems. Our approach leveraged parallel-in-time techniques based on \glspl{dft}, yielding a family of block-diagonalizable preconditioners that are well-suited for high-performance parallel implementations. This methodology allows the preconditioner to tend towards an exact approximation in the zero-viscosity limit.
We further presented efficient solvers for the systems arising at each time step, including recursive direct methods for structured grids with periodic and Neumann boundary conditions, and more general-purpose solvers for irregular grids.
Numerical experiments confirm the efficiency and scalability of the proposed preconditioners on the individual linear systems, highlighting their potential for large-scale optimal control simulations.

Future work may explore extensions of this approach to other discretizations in time, such as the Jordan--Kinderlehrer--Otto scheme (JKO) time-stepping scheme \cite{MR1617171}, which has received recent attention in the literature as it has been applied to solve nonlinear Fokker--Planck equations by using, for instance, primal--dual algorithms \cite{carrillo_primal_2022}. 
Other possible extensions of our work include \glspl{mfg} with non-local couplings \cite{ACFK17}, the \gls{mfg} planning problem (where a target density is imposed in the constraints \cite{MR2888257,MR3195848}), and \glspl{mfg} with congestion~\cite{MR3765549}. The solution method can be improved by including accelerators from the popular primal--dual approach for linear programs, extending the parallel implementation to usage on GPUs, or applying preconditioned iterative solvers within different outer iterative methods such as ADMM or Newton's method.

\bmhead{Acknowledgements}
HWL acknowledges financial support from the School of Mathematics at the University of Edinburgh, the Research Council Faroe Islands grant FRC0476, and the Betri Foundation. JWP acknowledges financial support from the Engineering and Physical Sciences Research Council (EPSRC) UK grants EP/S027785/1 and EP/Z533786/1. FJS was partially supported by l’Agence
Nationale de la Recherche (ANR), project ANR-22-CE40-0010, by KAUST through the subaward
agreement ORA-2021-CRG10-4674.6, and by Minist\`ere de l'Europe et des Affaires \'etrang\`eres (MEAE), project MATH AmSud 23-MATH-17. 
We thank Bernhard Heinzelreiter for useful discussions about parallelizable numerics. 
This work has made use of the resources provided by the Edinburgh Compute and Data Facility (ECDF).

\section*{Declarations}
\bmhead{Code availability} The code used for this paper is available on \href{https://github.com/heidiwol/PrecMFG}{GitHub}.
\bmhead{Data availability} This paper has no associated data.
\bmhead{Competing interests} The authors have no relevant financial or non-financial interests to disclose.
\bmhead{Funding} HWL was supported by the School of Mathematics at the University of Edinburgh, the Research Council Faroe Islands grant FRC0476, and the Betri Foundation. JWP was supported by the Engineering and Physical Sciences Research Council (EPSRC) UK grants EP/S027785/1 and EP/Z533786/1. FJS was partially supported by l’Agence
Nationale de la Recherche (ANR), project ANR-22-CE40-0010, by KAUST through the subaward
agreement ORA-2021-CRG10-4674.6, and by Minist\`ere de l'Europe et des Affaires \'etrang\`eres (MEAE), project MATH AmSud 23-MATH-17. 






\bibliography{biblio_MFG_CEMRACS}

\end{document}